\documentclass[final,onefignum,onetabnum]{siamart251216}

\usepackage{amssymb}
\usepackage{mathtools}
\usepackage{bm}
\usepackage{booktabs}
\usepackage{array}
\usepackage{float}
\usepackage{placeins}
\usepackage{graphicx}
\usepackage{hyperref}
\usepackage{enumitem}
\usepackage{microtype}
\usepackage{needspace}
\usepackage[normalem]{ulem}

\setlist[itemize]{leftmargin=2em}

\definecolor{red2}{RGB}{204,0,0}
\definecolor{blue2}{RGB}{0,103,165}
\definecolor{grey}{RGB}{100,100,100}
\newcommand{\RV}[1]{#1}
\hypersetup{colorlinks,linkcolor=blue2,citecolor=red2,urlcolor=blue2}
\newsiamremark{assumption}{Assumption}
\newsiamremark{remark}{Remark}
\newsiamremark{example}{Example}
\newsiamremark{question}{Question}
\newsiamremark{computation}{Computation}
\newsiamremark{result}{Result}
\numberwithin{equation}{section}

\newcommand{\R}{\mathbb{R}}

\newcommand{\diff}{\mathop{}\!\mathrm{d}}
\newcommand{\T}{{\top}}

\newcommand{\bPhi}{\boldsymbol{\Phi}}
\newcommand{\bPsi}{\boldsymbol{\Psi}}
\newcommand{\balpha}{\boldsymbol{\alpha}}

\newcommand{\bxi}{\boldsymbol{\xi}}

\newcommand{\bzeta}{\boldsymbol{\zeta}}
\newcommand{\bx}{\mathbf{x}}
\newcommand{\bq}{\mathbf{q}}
\newcommand{\bw}{\mathbf{w}}

\title{A Finite-Lag Markov Chain Model for Coarse-Graining Linear SDEs}
\author{Yuan~Chen\thanks{Division of Applied Mathematics, Brown University,
Providence, RI 02912 USA (\email{yuan\_chen1@brown.edu}).}}
\headers{Finite-Lag Markov Chain Model for Coarse-Graining}{Yuan Chen}

\begin{document}
\raggedbottom
\maketitle

\begin{abstract}
In this paper, we analyze a finite-lag Markov chain approximation for
coarse-grained linear stochastic differential equations. The chain
uses the limiting conditional transition law of the resolved variables to
advance on a prescribed time grid. Its stability and statistical properties
are studied, with particular attention to the role of the lag time. The chain preserves the exact equilibrium distribution and matches the autocovariance at the chosen lag, while errors at subsequent times are quantified. For a slow--fast system, the analysis relates lag selection to scale separation and dynamical accuracy. A continuous embedding of the chain enables a unified comparison with some existing continuous-time Markovian effective models based on the temporal statistics they preserve. Numerical examples illustrate the theoretical results.
\end{abstract}



\section{Introduction}
Coarse-graining seeks to describe the dynamics of a complex system using a small set of resolved variables. Such models are widely used in molecular and biomolecular simulation, including
protein folding and molecular kinetics, as well as in modeling low-frequency
climate variability \cite{noid2013perspective,best2010coordinate,prinz2011markov,majda2003systematic}.

The resolved variables generally do not form a Markov process, even when the
full system is Markovian. The Mori--Zwanzig formalism represents their evolution
through memory terms and fluctuating forces
\cite{mori1965transport,zwanzig1973nonlinear,chorin2000optimal,lin2021data}.
This representation underlies non-Markovian effective models and numerical
procedures for computing their memory kernels
\cite{darve2009computing,hijon2010mori,li2015incorporation,hudson2020coarse}.
Such non-Markovian descriptions include generalized Langevin models, in which the memory kernel and the fluctuating force are linked by fluctuation–dissipation relations \cite{hijon2010mori,klippenstein2021introducing}.
Methods for estimating these models from microscopic simulations include
kernel parametrization, reconstruction from correlation functions, and
auxiliary-variable representations
\cite{lei2016data,jung2018generalized,grogan2020data}.
Nonlinear models with position- or state-dependent memory have also been
derived and learned from data \cite{vroylandt2022position,ge2024data}.
Klippenstein et al.\ \cite{klippenstein2021introducing} review memory effects
and their treatment in coarse-grained molecular simulations.

Alternatively, Markovian approaches approximate the resolved dynamics by a Markov process on the resolved state space. Such approximations are often appropriate when the unresolved variables relax much faster than the resolved dynamics. Averaging and homogenization provide effective equations in this regime \cite{majda2001mathematical,pavliotis2008multiscale}, with related
computational frameworks for estimating coarse dynamics from microscopic
simulations \cite{kevrekidis2003equation,e2007heterogeneous}.
Conditional-expectation constructions and their approximation errors have
been studied for diffusions, including
pathwise estimates for nonlinear reaction coordinates
\cite{legoll2010effective,legoll2017pathwise,lelievre2019pathwise,legoll2019effective}.
Other constructions fit coarse dynamics through path-space relative entropy
\cite{harmandaris2016path}, or learn drift, diffusion, and friction coefficients
from data \cite{boninsegna2018sparse,sachs2025machine}.
The dynamical properties of effective models have also been analyzed through
dominant timescales, spectral approximation, and transition rates
\cite{zhang2016effective,nuske2021spectral}.
Closely related to the present work, Hudson and Li
\cite{hudson2024dynamical} analyze the dynamical statistics of coarse-grained
linear SDEs and introduce a memory-corrected Markovian approximation.
Alongside continuous-time effective equations, transition statistics over a prescribed lag time provide another way to characterize the resolved dynamics, see e.g. Markov state models
\cite{schutte1999direct,prinz2011markov}, time-lagged independent component
analysis (TICA) \cite{perez2013identification}, and the variational approach
to conformational dynamics (VAC) \cite{noe2013variational,klus2018data}.
A related construction defines a Markov chain on the resolved state space
using the equilibrium conditional transition law at a prescribed lag.
Its operator formulation was studied by Chekroun et al.\
\cite{chekroun2020ruelle}, while Zhang and Sch{\"u}tte
\cite{zhang2025finding} analyze the resulting effective dynamics for general
collective variables. From a data-driven perspective, we developed a finite-lag Markov model for
multiscale stochastic systems by learning conditional transitions from
resolved data \cite{chen2024datadriven,chen2024learning}.

\RV{In this paper, we analyze this finite-lag Markov construction for the
resolved dynamics of linear SDEs. For a prescribed lag time $\Delta>0$,
the chain advances on the grid $0,\Delta,2\Delta,\ldots$ by sampling at each
step from the long-time limit of the lag-$\Delta$ conditional transition law,
conditioned on the current resolved state. We study its stability and
statistical properties, quantify dynamical errors, and examine the dependence
on the lag time and scale separation. A continuous-time embedding relates
the chain to existing continuous-time Markovian effective models. The main
results are as follows.}
\begin{itemize}[topsep=0pt,partopsep=0pt]
    \item \RV{\textit{Long-time statistics.} We prove stability and
    convergence to the exact resolved equilibrium distribution for every
    $\Delta>0$. The chain also reproduces the equilibrium autocovariance at
    lag $\Delta$. At longer lags, we quantify the discrepancy introduced by
    repeated transitions and show that the chain underestimates the
    discrete integrated autocovariance, with an explicit error bound. For the
    two-dimensional slow--fast system, the relative autocovariance error
    converges to zero as the scale separation tends to infinity, on any fixed finite horizon with fixed $\Delta>0$.}

    \item \RV{\textit{Role of the lag time.} We characterize how $\Delta$
    affects the predicted temporal correlations, including a monotonic
    decrease in the eigenvalues of the normalized autocovariance under grid
    refinement at a fixed time. For a two-dimensional slow--fast system, we
    determine how the error depends on the lag time relative to the fast
    relaxation time and derive a sufficient condition on $\Delta$ to meet
    a prescribed error tolerance over a finite horizon.}

    \item \RV{\textit{Comparison with continuous-time Markovian models.} We construct a
    continuous embedding that reproduces the law of the chain on the grid.
    This provides a unified framework for comparing the finite-lag chain with
    generator-based and memory-corrected continuous-time effective models
    through the local or global properties of the normalized autocovariance they
    preserve. Under the framework, we prove convergence of the
    embedding to the generator-based model as $\Delta\to0$. For a scalar
    resolved variable, we identify lag times at which the chain and the
    memory-corrected model agree in law on the grid.}
\end{itemize}

The remainder of the paper is organized as follows.
Section~\ref{sec:setup} introduces the setting and the finite-lag Markov chain.
Section~\ref{sec:properties} establishes its statistical properties, and
Section~\ref{sec:2d} analyzes a two-dimensional slow-fast system
and the choice of lag time. Section~\ref{sec:comparison} relates the chain
to existing continuous-time Markovian effective models. Section~\ref{sec:examples} presents
numerical examples, and Section~\ref{sec:conclu} concludes the paper.

\Needspace{6\baselineskip}
\section{Setup and Objective}\label{sec:setup}
\subsection{Setup}
We consider the linear SDE on $\R^N$ with
potential function $V(\bq)=\frac12\bq^\top\mathbf{K}\bq$:
\begin{equation}\label{model}
  \diff\bq_t=-\mathbf{K}\bq_t\,\diff t+\sqrt{2\beta^{-1}}\,\diff\bw_t,
\end{equation}
where $N$ is the dimension of the full state space, $\beta>0$ is the inverse temperature, and assume $\mathbf{K}\in\R^{N\times N}$ is symmetric and positive definite, $\bw_t$ is an $N$-dimensional standard Brownian motion.

Let $n<N$ be the dimension of resolved variables, and set $m=N-n$. We let $\widetilde{\bPhi}\in\R^{n\times N}$ be a full-rank linear coarse-graining
matrix. We define the normalized coarse-graining matrix and the resolved
variable by
\begin{equation}\label{coarsemap}
  \bPhi
  :=
  \left(\widetilde{\bPhi}\widetilde{\bPhi}^{\top}\right)^{-1/2}\widetilde{\bPhi},
  \qquad
  \bxi_t:=\bPhi\bq_t .
\end{equation}
Then $\bPhi\bPhi^\top=\mathbf{I}_n$, where $\mathbf{I}_n$ is the identity matrix with dimension $n$. Choose
$\bPsi\in\R^{m\times N}$ with orthonormal rows spanning the orthogonal
complement of the row space of $\bPhi$. Thus
\begin{equation}\label{complement}
  \bPsi\bPsi^\top=\mathbf{I}_m,
  \qquad
  \bPhi\bPsi^\top=0,
  \qquad
  \mathbf{I}_N=\bPhi^\top\bPhi+\bPsi^\top\bPsi .
\end{equation}
Set $\bzeta_t:=\bPsi\bq_t$, $\diff\bw_{0,t}:=\bPhi \diff\bw_t$, and
$\diff\bw_{1,t}:=\bPsi \diff\bw_t$. By construction, $\bw_0$ and $\bw_1$ are independent
standard Brownian motions. In the coordinates $\bx=(\bxi,\bzeta)$, the SDE can
therefore be written as
\begin{equation}
  \begin{pmatrix} \diff\bxi_t \\ \diff\bzeta_t \end{pmatrix}
  = -\begin{pmatrix} \mathbf{A}_0 & \balpha \\ \balpha^\T & \mathbf{A}_1 \end{pmatrix}
    \begin{pmatrix} \bxi_t \\ \bzeta_t \end{pmatrix}\diff t
  + \sqrt{2\beta^{-1}}\begin{pmatrix} \diff\bw_{0,t} \\ \diff\bw_{1,t} \end{pmatrix},
  \label{model_block}
\end{equation}
where
$\mathbf{A}_0=\bPhi\mathbf{K}\bPhi^\top\in \R^{n\times n}$,
$\mathbf{A}_1=\bPsi\mathbf{K}\bPsi^\top\in \R^{m\times m}$, and
$\balpha=\bPhi\mathbf{K}\bPsi^\top\in \R^{n\times m}$. We denote the full block
matrix in \eqref{model_block} by
\begin{equation}\label{A}
    \mathbf{A}
    :=
    \begin{pmatrix}
        \mathbf{A}_0 & \balpha \\
        \balpha^\top & \mathbf{A}_1
    \end{pmatrix}.
\end{equation}
The initial condition is Gaussian,
$(\bxi_0,\bzeta_0)\sim \mathcal{N}(\boldsymbol{\mu}_0,\boldsymbol{\Sigma}_0)$,
where $\boldsymbol{\mu}_0\in\R^N$ is the initial mean vector and
$\boldsymbol{\Sigma}_0\in\R^{N\times N}$ is the initial covariance matrix. With
respect to the block decomposition $\bx=(\bxi,\bzeta)$, we write
$[\boldsymbol{\mu}_0]_{\xi}$ for the $\bxi$-component of
$\boldsymbol{\mu}_0$, and $[\boldsymbol{\Sigma}_0]_{\xi\xi}$ for the
corresponding $\bxi$-$\bxi$ covariance block. We denote the equilibrium
distribution by $\pi_{\mathbf{x}}$ and its $\xi$-marginal by $\pi_\xi$.

We usually have $n \ll N$, which means the resolved variables $\bxi_t$ are much fewer than the unresolved variables $\bzeta_t$. We study a finite-lag Markov approximation of the resolved dynamics and analyze its long-time statistical properties.

\subsection{Time-Dependent Conditional Transition Law}\label{sec:transition}
The resolved variables are not expected to form a closed Markovian system. Although the full process $\mathbf{x}_t=(\bxi_t,\bzeta_t)$ is Markovian, the value of $\bxi_t$ alone does not determine how the unresolved variable $\bzeta_t$ is distributed. Consequently, the lag-$\Delta$ conditional transition law of $\bxi_{t+\Delta}$ given $\bxi_t$ still depends on the conditional law of $\bzeta_t$ at the same time. Since this conditional law evolves with the distribution of the full process, the resolved lag-$\Delta$ conditional transition law is generally time-dependent.

In the present linear Gaussian setting, we can locate this time-dependence explicitly. We first record the Schur complements
\begin{equation}\label{Schur}
    \mathbf{B}_0 = \mathbf{A}_0 - \boldsymbol{\alpha}\mathbf{A}_1^{-1}\boldsymbol{\alpha}^\top,
    \qquad
    \mathbf{B}_1 = \mathbf{A}_1 - \boldsymbol{\alpha}^\top\mathbf{A}_0^{-1}\boldsymbol{\alpha},
\end{equation}
both of which are symmetric positive definite. The block inverse of $\mathbf{A}$ is
\begin{equation}\label{Ainv}
    \mathbf{A}^{-1}
    =
    \begin{pmatrix}
        \mathbf{B}_0^{-1}
        &
        -\mathbf{A}_0^{-1}\boldsymbol{\alpha}\mathbf{B}_1^{-1}
        \\
        -\mathbf{B}_1^{-1}\boldsymbol{\alpha}^\top\mathbf{A}_0^{-1}
        &
        \mathbf{B}_1^{-1}
    \end{pmatrix}.
\end{equation}
The Gaussian structure is preserved by the linear dynamics. Since $\mathbf{A}$ is symmetric positive definite and the initial condition is Gaussian, $\mathbf{x}_t$ has mean and covariance
\begin{equation} \label{moments}
    \boldsymbol{\mu}(t) = e^{-\mathbf{A}t}\boldsymbol{\mu}_0,
    \qquad
    \boldsymbol{\Sigma}(t)
    =
    e^{-\mathbf{A}t}\!\left(\boldsymbol{\Sigma}_0 - \beta^{-1}\mathbf{A}^{-1}\right)e^{-\mathbf{A}t}
    +
    \beta^{-1}\mathbf{A}^{-1}.
\end{equation}
We decompose these quantities according to $\mathbf{x}=(\bxi,\bzeta)$:
\begin{equation}\label{moments_block}
    \boldsymbol{\mu}(t) = \begin{pmatrix} \boldsymbol{\mu}_\xi(t) \\ \boldsymbol{\mu}_\zeta(t) \end{pmatrix},
    \qquad
    \boldsymbol{\Sigma}(t)
    =
    \begin{pmatrix}
        \boldsymbol{\Sigma}_{\xi\xi}(t) & \boldsymbol{\Sigma}_{\xi\zeta}(t) \\
        \boldsymbol{\Sigma}_{\zeta\xi}(t) & \boldsymbol{\Sigma}_{\zeta\zeta}(t)
    \end{pmatrix}.
\end{equation}
As $t\to\infty$,
$\boldsymbol{\mu}(t) \to \mathbf{0}$ and $ \boldsymbol{\Sigma}(t) \to \beta^{-1}\mathbf{A}^{-1}$.

Over a lag time $\Delta$, the full dynamics satisfies
\begin{equation} \label{trans_full}
    \mathbf{x}_{t+\Delta}
    =
    e^{-\mathbf{A}\Delta}\mathbf{x}_t
    +
    \boldsymbol{\eta}_\Delta,
    \qquad
    \boldsymbol{\eta}_\Delta \sim \mathcal{N}(\mathbf{0}, \sigma_\Delta),
\end{equation}
where $\boldsymbol{\eta}_\Delta$ is independent of $\mathbf{x}_t$ and
$\sigma_\Delta = \beta^{-1}\mathbf{A}^{-1}(\mathbf{I}_N - e^{-2\mathbf{A}\Delta})$.
With the block notation
\begin{equation} \label{expA}
    e^{-\mathbf{A}\Delta}
    =
    \begin{pmatrix}
        \mathbf{E}_0 & \mathbf{e} \\
        \mathbf{e}^\top & \mathbf{E}_1
    \end{pmatrix}.
\end{equation}
The resolved component of \eqref{trans_full} is
\begin{equation}\label{trans_resolved}
    \boldsymbol{\xi}_{t+\Delta}
    =
    \mathbf{E}_0\boldsymbol{\xi}_t
    +
    \mathbf{e}\boldsymbol{\zeta}_t
    +
    \boldsymbol{\eta}_\xi,
\end{equation}
where $\boldsymbol{\eta}_\xi$ is the $\boldsymbol{\xi}$-component of $\boldsymbol{\eta}_\Delta$. Thus the resolved transition is obtained by averaging the unresolved variable in \eqref{trans_resolved} against its conditional law given $\boldsymbol{\xi}_t$. Since the joint law is Gaussian, this gives the Gaussian transition
\begin{align}\label{moments_cond}
    \mathbb{E}[\boldsymbol{\xi}_{t+\Delta}\mid\boldsymbol{\xi}_t]
    &= \mathbf{E}_0\boldsymbol{\xi}_t + \mathbf{e}\,\mathbb{E}[\boldsymbol{\zeta}_t\mid\boldsymbol{\xi}_t], \\
    \mathrm{Cov}(\boldsymbol{\xi}_{t+\Delta}\mid\boldsymbol{\xi}_t)
    &= [\sigma_\Delta]_{\xi\xi}
    +
    \mathbf{e}\,\mathrm{Cov}(\boldsymbol{\zeta}_t\mid\boldsymbol{\xi}_t)\,\mathbf{e}^\top,
\end{align}
where, for $t>0$, the conditional moments are
\begin{align}\label{mean_cond}
    \mathbb{E}[\boldsymbol{\zeta}_t\mid\boldsymbol{\xi}_t]
    &=
    \boldsymbol{\mu}_\zeta(t)
    +
    \boldsymbol{\Sigma}_{\zeta\xi}(t)\boldsymbol{\Sigma}_{\xi\xi}(t)^{-1}
    \bigl(\boldsymbol{\xi}_t-\boldsymbol{\mu}_\xi(t)\bigr), \\
    \mathrm{Cov}(\boldsymbol{\zeta}_t\mid\boldsymbol{\xi}_t)
    &=
    \boldsymbol{\Sigma}_{\zeta\zeta}(t)
    -
    \boldsymbol{\Sigma}_{\zeta\xi}(t)\boldsymbol{\Sigma}_{\xi\xi}(t)^{-1}\boldsymbol{\Sigma}_{\xi\zeta}(t).\label{cov_cond}
\end{align}
Equations \eqref{moments_cond}--\eqref{mean_cond} make this dependence explicit in the present model: the exact lag-$\Delta$ conditional transition law involves the time-dependent quantities $\boldsymbol{\mu}(t)$ and $\boldsymbol{\Sigma}(t)$. Thus $\mathbb{P}(\boldsymbol{\xi}_{t+\Delta}\mid\boldsymbol{\xi}_t)$ generally depends on $t$, which motivates a time-homogeneous approximation.

\subsection{Finite-Lag Markov Chain Approximation} \label{sec:markov}
The time-dependence displayed above can be handled in different ways. One approach is to retain the effect of the unresolved variables through memory terms, leading to generalized Langevin equations (GLEs) \cite{zwanzig1973nonlinear,zwanzig2001nonequilibrium}. In such models, the influence of $\bzeta_t$ is represented through history-dependent drift and noise terms on the resolved variables. Another approach is to replace the time-dependent lag-$\Delta$ conditional transition law by a time-homogeneous one-step transition law on the resolved variables alone. Such Markov chain approximations are most useful when the unresolved variables relax on a faster time scale \cite{hudson2024dynamical,legoll2019effective,legoll2010effective}. In this work, we focus on the latter approach.

For a fixed lag time $\Delta>0$, we approximate the resolved dynamics by a discrete-time Markov chain $\{\widetilde{\bxi}_{k\Delta}\}_{k\geq 0}$ on $\R^n$. Its initial distribution is chosen to match that of the resolved variable,
$\widetilde{\bxi}_0 \sim \mathcal{N}\bigl([\boldsymbol{\mu}_0]_{\xi},[\boldsymbol{\Sigma}_0]_{\xi\xi}\bigr)$.
The one-step transition law of the chain is chosen to match the limiting lag-$\Delta$ conditional transition law computed in Section~\ref{sec:transition}:
\begin{equation}\label{kernel}
  \mathbb{P}_{\infty}(\widetilde{\bxi},B;\Delta)
  :=
  \lim_{s \to \infty}
  \mathbb{P}\bigl(\boldsymbol{\xi}_{s+\Delta}\in B\mid \boldsymbol{\xi}_s=\widetilde{\bxi}\,\bigr),
  \qquad B\subset\R^n\ \text{Borel}.
\end{equation}
Unless otherwise specified, $\Delta$ is fixed throughout the analysis. The quantities derived from this transition law depend on $\Delta$, but we omit this dependence from the notation when no comparison across different lag times is being made.
The limit exists because $\mathbf{A}$ is symmetric positive definite. Using the limiting mean and covariance above, \eqref{moments_cond}--\eqref{cov_cond}, and the block inverse formula \eqref{Ainv}, this limiting transition law is Gaussian with
\begin{equation} \label{trans_limit}
    \mathbb{E}_{\infty}[\widetilde{\boldsymbol{\xi}}_{t+\Delta}\mid\widetilde{\boldsymbol{\xi}}_t\,]
    =
    \left(\mathbf{E}_0-\mathbf{e}\mathbf{A}_1^{-1}\boldsymbol{\alpha}^\top\right)\widetilde{\boldsymbol{\xi}}_t,
    \qquad
    \mathrm{Cov}_{\infty}(\widetilde{\boldsymbol{\xi}}_{t+\Delta}\mid\widetilde{\boldsymbol{\xi}}_t\,)
    =
    [\sigma_\Delta]_{\xi\xi}
    +
    \mathbf{e}\,\beta^{-1}\mathbf{A}_1^{-1}\,\mathbf{e}^\top .
\end{equation}
We denote the one-step conditional-mean coefficient and one-step conditional covariance in \eqref{trans_limit} by
\begin{equation}\label{FG}
    \mathbf{F}_\infty
    =
    \mathbf{E}_0-\mathbf{e}\mathbf{A}_1^{-1}\boldsymbol{\alpha}^\top \in \mathbb{R}^{n\times n},
    \qquad
    \mathbf{G}_\infty
    =
    [\sigma_\Delta]_{\xi\xi}
    +
    \beta^{-1}\mathbf{e}\mathbf{A}_1^{-1}\mathbf{e}^\top \in \mathbb{R}^{n\times n}.
\end{equation}
We then define the finite-lag Markov chain approximation by setting, for every Borel set $B\subset\R^n$,
\begin{equation}\label{transMC}
    \mathbb{P}\bigl(\widetilde{\bxi}_{(k+1)\Delta}\in B
    \mid \widetilde{\bxi}_{k\Delta}=\widetilde{\bxi}\,\bigr)
    =
    \mathbb{P}_{\infty}(\widetilde{\bxi},B),
    \qquad k\geq0.
\end{equation}
Following this notation, the same Markov chain can be equivalently defined by
the recursion
\begin{equation}\label{MC}
    \widetilde{\boldsymbol{\xi}}_{j\Delta}
    =
    \mathbf{F}_\infty\widetilde{\boldsymbol{\xi}}_{(j-1)\Delta}
    +
    \boldsymbol{\eta}_j,
    \qquad
    \boldsymbol{\eta}_j\sim\mathcal{N}(\mathbf{0},\mathbf{G}_\infty),
    \quad j\ge1,
\end{equation}
where $\boldsymbol{\eta}_1,\boldsymbol{\eta}_2,\ldots$ are i.i.d.\ and
independent of $\widetilde{\boldsymbol{\xi}}_0$.

\section{Statistical Properties}\label{sec:properties}
In this section, we analyze the statistical properties of the finite-lag Markov chain approximation.
\subsection{Stability and Equilibrium Properties}\label{sec:equilibrium}
We first identify the equilibrium distribution that the effective dynamics should
preserve. Since $\mathbf{A}$ is symmetric positive definite, the full linear
dynamics has the unique invariant Gaussian distribution
$\pi_{\mathbf{x}} = \mathcal{N}(\mathbf{0},\beta^{-1}\mathbf{A}^{-1})$.
Using the block inverse formula \eqref{Ainv}, the corresponding marginal
equilibrium distribution of the resolved variable is
$\pi_\xi = \mathcal{N}(\mathbf{0},\beta^{-1}\mathbf{B}_0^{-1})$.
Thus a basic consistency requirement for the finite-lag Markov chain is that it
has $\pi_\xi$ as its equilibrium distribution.
We prove the stronger statement that the transition law converges to this
distribution from every fixed initial state. 

Using the recursive form \eqref{MC}, induction gives
\begin{equation}\label{MC_k}
    \widetilde{\boldsymbol{\xi}}_{k\Delta}
    =
    \mathbf{F}_\infty^k\widetilde{\boldsymbol{\xi}}_0
    +
    \sum_{j=0}^{k-1}\mathbf{F}_\infty^j\boldsymbol{\eta}_{k-j}.
\end{equation}
Therefore the $k$-step transition law is
\begin{equation}\label{transMC_k}
    \mathbb{P}(\widetilde{\boldsymbol{\xi}}_{k\Delta}|\widetilde{\boldsymbol{\xi}}_0\,)
    =
    \mathcal{N}\!\left(
    \mathbf{F}_\infty^k\widetilde{\boldsymbol{\xi}}_0,
    \sum_{j=0}^{k-1}
    \mathbf{F}_\infty^j\mathbf{G}_\infty(\mathbf{F}_\infty^\top)^j
    \right).
\end{equation}
The next lemma, proved in Appendix~\ref{app:statistical}, verifies this stability for every fixed $\Delta>0$. Throughout the proof and the subsequent discussion, we use the Loewner order to compare symmetric matrices: $\mathbf{M}\prec\mathbf{N}$
means that $\mathbf{N}-\mathbf{M}$ is positive definite, and
$\mathbf{M}\preceq\mathbf{N}$ means that $\mathbf{N}-\mathbf{M}$ is positive
semidefinite.

\begin{lemma}[Fixed-lag stability]\label{lem:stability}
For every fixed $\Delta>0$,
\begin{equation}\label{stability}
    \rho(\mathbf{F}_\infty)<1 .
\end{equation}
\end{lemma}

Lemma~\ref{lem:stability} shows that the finite-lag Markov chain approximation is stable for every fixed lag time $\Delta>0$. The next theorem shows that the transition law converges to the equilibrium distribution from every initial state.

\begin{theorem}[Convergence to equilibrium]\label{thm:equilibrium}
For every fixed $\Delta>0$ and every initial state
$\widetilde{\boldsymbol{\xi}}\in\mathbb{R}^n$, the transition law of
$\widetilde{\boldsymbol{\xi}}_{k\Delta}$ conditional on
$\widetilde{\boldsymbol{\xi}}_0=\widetilde{\boldsymbol{\xi}}$ converges weakly
to $\pi_\xi$ as $k\to\infty$.
\end{theorem}

\begin{proof}
Fix $\widetilde{\boldsymbol{\xi}}_0=\widetilde{\boldsymbol{\xi}}$.
By \eqref{transMC_k}, the $k$-step transition law is Gaussian with mean
$\mathbf{m}_k = \mathbf{F}_\infty^k\widetilde{\boldsymbol{\xi}}$,
and covariance
\begin{equation}\label{covMC}
    \boldsymbol{\Sigma}_k
    =
    \sum_{j=0}^{k-1}
    \mathbf{F}_\infty^j\mathbf{G}_\infty(\mathbf{F}_\infty^\top)^j .
\end{equation}
Since $\rho(\mathbf{F}_\infty)<1$ by
Lemma~\ref{lem:stability}, $\mathbf{F}_\infty^k\to\mathbf{0}$ as
$k\to\infty$, and therefore $\mathbf{m}_k\to\mathbf{0}$ as $k\to\infty$. For
the covariance, the identity established in Appendix~\ref{app:statistical}
gives
\begin{equation}\label{Gidentity}
    \mathbf{G}_\infty
    =
    \beta^{-1}\mathbf{B}_0^{-1}
    -
    \mathbf{F}_\infty\beta^{-1}\mathbf{B}_0^{-1}\mathbf{F}_\infty^\top .
\end{equation}
Substituting \eqref{Gidentity} into \eqref{covMC}, and using a telescoping sum, we have
\begin{align}\label{covMC_limit}
    \boldsymbol{\Sigma}_k
    &=
    \sum_{j=0}^{k-1}
    \mathbf{F}_\infty^j
    \left(
    \beta^{-1}\mathbf{B}_0^{-1}
    -
    \mathbf{F}_\infty\beta^{-1}\mathbf{B}_0^{-1}\mathbf{F}_\infty^\top
    \right)
    (\mathbf{F}_\infty^\top)^j \notag\\
    &=
    \beta^{-1}\mathbf{B}_0^{-1}
    -
    \mathbf{F}_\infty^k
    \beta^{-1}\mathbf{B}_0^{-1}
    (\mathbf{F}_\infty^\top)^k
    \to
    \beta^{-1}\mathbf{B}_0^{-1},
    \qquad k\to\infty.
\end{align}
The Gaussian transition laws are determined by their means and covariances, so
they converge weakly to
$\mathcal{N}(\mathbf{0},\beta^{-1}\mathbf{B}_0^{-1})=\pi_\xi$ as
$k\to\infty$. In particular, $\pi_\xi$ is the unique equilibrium distribution
of the Markov chain.
\end{proof}

\subsection{Autocovariance Function}\label{sec:acf}
At equilibrium,
$\mathbb{E}_{\pi_{\mathbf x}}[\bx_0\bx_0^\top]
=\beta^{-1}\mathbf{A}^{-1}$. Using the exact solution of the full dynamics,
the autocovariance of the resolved variables is
\begin{align}\label{R}
    \mathbf{R}(t)
    :=
    \mathbb{E}_{\pi_{\mathbf x}}
    \left[
    \boldsymbol{\xi}_t\boldsymbol{\xi}_0^\top
    \right]
    =
    \left[
    e^{-\mathbf{A}t}
    \mathbb{E}_{\pi_{\mathbf x}}[\bx_0\bx_0^\top]
    \right]_{\xi\xi}=
    \beta^{-1}
    \left[e^{-\mathbf{A}t}\mathbf{A}^{-1}\right]_{\xi\xi}.
\end{align}
\RV{Since
$\mathbf{R}(0)=\beta^{-1}\mathbf{B}_0^{-1}$, we normalize the
autocovariance by its equilibrium value. This maps the equilibrium covariance
to $\mathbf{I}_n$ and, through the compression formula below, makes the
symmetric contraction structure available for the matrix analysis.} We define
the normalized equilibrium autocovariance by
\begin{equation}\label{H}
    \mathbf{H}(t)
    :=
    \beta\mathbf{B}_0^{1/2}\mathbf{R}(t)\mathbf{B}_0^{1/2},
    \qquad
    \mathbf{R}(t)
    =
    \beta^{-1}\mathbf{B}_0^{-1/2}
    \mathbf{H}(t)\mathbf{B}_0^{-1/2}.
\end{equation}
\RV{This normalization will also place the models compared in
Section~\ref{sec:comparison} in a common coordinate system.} We set
\begin{equation}\label{U}
    \mathbf{U}
    :=
    \mathbf{A}^{-1/2}
    \begin{pmatrix}
        \mathbf{I}_n\\
        \mathbf{0}_{m\times n}
    \end{pmatrix}
    \mathbf{B}_0^{1/2},
    \qquad
    \mathbf{U}^\top\mathbf{U}=\mathbf{I}_n,
    \qquad
    \mathbf{H}(t)
    =
    \mathbf{U}^\top e^{-\mathbf{A}t}\mathbf{U}.
\end{equation}
\RV{For every $\Delta>0$, \eqref{H} gives}
$\mathbf{H}(\Delta) = \mathbf{B}_0^{1/2}\mathbf{F}_\infty\mathbf{B}_0^{-1/2}$.
\RV{When the Markov chain is initialized at its equilibrium distribution
$\pi_\xi$, stationarity gives, for every $k\geq0$,}
\begin{equation}\label{Rmc}
    \widetilde{\mathbf{R}}(k\Delta)
    =
    \mathbb{E}\!
    \left[
    \widetilde{\boldsymbol{\xi}}_{(j+k)\Delta}
    \widetilde{\boldsymbol{\xi}}_{j\Delta}^\top
    \right] 
    =
    \mathbf{F}_\infty^k\,
    \mathbb{E}
    \left[
    \widetilde{\boldsymbol{\xi}}_{j\Delta}
    \widetilde{\boldsymbol{\xi}}_{j\Delta}^\top
    \right]=\mathbf{F}_\infty^k\beta^{-1}\mathbf{B}_0^{-1},
\end{equation}
\RV{The last equality uses
$\mathbb{E}_{\pi_\xi}[\widetilde{\boldsymbol\xi}_{j\Delta}
\widetilde{\boldsymbol\xi}_{j\Delta}^{\top}]
=\beta^{-1}\mathbf B_0^{-1}$.}
\RV{Consequently, its normalized autocovariance satisfies the iterative
identity}
\begin{equation}\label{Hmc}
    \widetilde{\mathbf{H}}(k\Delta)
    :=
    \beta\mathbf{B}_0^{1/2}
    \widetilde{\mathbf{R}}(k\Delta)
    \mathbf{B}_0^{1/2}
    =
    \left(
    \mathbf{B}_0^{1/2}\mathbf{F}_\infty\mathbf{B}_0^{-1/2}
    \right)^k
    =
    \mathbf{H}(\Delta)^k.
\end{equation}

\subsubsection{Autocovariance Error}
\RV{The effective chain matches the exact autocovariance after one transition.
At later grid points, its error is generated by repeatedly replacing the
resolved dynamics by this one-step transition law. We first record the exact
one-step property and then quantify the accumulated error.}
\begin{proposition}[One-step exactness]\label{prop:one_step}
For every $\Delta>0$,
\begin{equation}\label{Rmatch}
    \widetilde{\mathbf{R}}(\Delta)=\mathbf{R}(\Delta).
\end{equation}
\end{proposition}

\begin{proof}
\RV{Equations \eqref{H} and
\eqref{Hmc} give
$\widetilde{\mathbf H}(\Delta)=\mathbf H(\Delta)$. Transforming back by
\eqref{H} proves \eqref{Rmatch}.}
\end{proof}

\RV{This exact finite-lag match will reappear in
Section~\ref{sec:comparison}. To measure the error accumulated over subsequent
grid points, define}
\begin{align}\label{int_error}
    \boldsymbol{\mathcal E}_\Delta
    &:=
    \Delta\sum_{k=0}^{\infty}
    \left[
    \mathbf{R}(k\Delta)
    -
    \widetilde{\mathbf{R}}(k\Delta)
    \right] =
    \Delta\sum_{k=2}^{\infty}
    \left[
    \mathbf{R}(k\Delta)
    -
    \widetilde{\mathbf{R}}(k\Delta)
    \right],
\end{align}
\RV{where the second equality follows from the agreement at lags $0$ and
$\Delta$. This series is absolutely convergent because
\[
    \|\mathbf{R}(k\Delta)\|_2
    +\|\widetilde{\mathbf{R}}(k\Delta)\|_2
    \le
    \beta^{-1}
    \left(
    \|\mathbf{A}^{-1}\|_2+\|\mathbf{B}_0^{-1}\|_2
    \right)
    e^{-\lambda_{\min}(\mathbf{A})k\Delta},
\]
where $\lambda_{\min}(\mathbf A)$ denotes the smallest eigenvalue of
$\mathbf A$. We also used that $\mathbf U$ is an isometry and hence
$\|\mathbf H(\Delta)\|_2\le e^{-\lambda_{\min}(\mathbf A)\Delta}$.}

\RV{Choose $\mathbf{V}\in\mathbb{R}^{N\times(N-n)}$ such that
$(\mathbf{U},\mathbf{V}) \in \mathbb{R}^{N\times N}$ is an orthogonal matrix. Theorem~\ref{thm:int_error} and Corollary~\ref{cor:int_bound} characterize
the sign and magnitude of the accumulated error, with proofs deferred to Appendix~\ref{app:statistical}.}

\begin{theorem}[Integrated autocovariance error]
\label{thm:int_error}
For every $\Delta>0$,
\begin{equation}\label{int_error_nonnegative}
    \boldsymbol{\mathcal E}_\Delta\succeq\mathbf{0}.
\end{equation}
\end{theorem}

\RV{The accumulated error in \eqref{int_error} reflects the failure of the
compressed dynamics to form a semigroup. Using \eqref{U}, its two-step
composition defect is}
\begin{equation}\label{eta}
    \RV{\eta_\Delta^2
    :=
    \left\|
    \mathbf{H}(2\Delta)-\mathbf{H}(\Delta)^2
    \right\|_2
    =
    \|\mathbf{V}^\top e^{-\mathbf{A}\Delta}\mathbf{U}\|_2^2.}
\end{equation}
\RV{To state the resulting error bound, also define}
\begin{equation}\label{bound_const}
\RV{\begin{aligned}
    \ell_\Delta:=
    \|\mathbf{V}^\top e^{-\mathbf{A}\Delta}\mathbf{V}\|_2,
    ~~~~
    h_\Delta:=\|\mathbf{H}(\Delta)\|_2=\|\mathbf{U}^\top e^{-\mathbf{A}\Delta}\mathbf{U}\|_2,
    ~~~~
    r_\Delta:=\|e^{-\mathbf{A}\Delta}\|_2,
\end{aligned}}
\end{equation}
which satisfy $0\le \max\{\ell_\Delta,h_\Delta\} \leq r_\Delta<1$ for every $\Delta>0$.

\begin{corollary}[Global error bound]\label{cor:int_bound}
For every $\Delta>0$,
\begin{equation}\label{int_bound}
    \|\boldsymbol{\mathcal E}_\Delta\|_2
    \le
    \beta^{-1}\Delta\|\mathbf{B}_0^{-1}\|_2
    \frac{\eta_\Delta^2}
    {(1-r_\Delta)(1-\ell_\Delta)(1-h_\Delta)}.
\end{equation}
In particular,
\begin{equation}\label{int_decay}
    \|\boldsymbol{\mathcal E}_\Delta\|_2
    \le
    \beta^{-1}\|\mathbf{B}_0^{-1}\|_2
    \frac{
    \Delta e^{-2\lambda_{\min}(\mathbf{A})\Delta}
    }{
    \left(1-e^{-\lambda_{\min}(\mathbf{A})\Delta}\right)^3
    },
\end{equation}
and hence $\boldsymbol{\mathcal E}_\Delta\to\mathbf{0}$ as
$\Delta\to\infty$.
\end{corollary}

\RV{Thus $\eta_\Delta^2$ measures the resolved--unresolved transfer that
breaks the semigroup property, while $(1-\ell_\Delta)^{-1}$ measures the
persistence of this information in the complementary subspace. The remaining
resolvent factors accumulate the defect under the full and effective dynamics.
The large-$\Delta$ limit above concerns the error on the increasingly sparse
grid $\{k\Delta\}_{k\ge0}$. It does not assert uniform accuracy between
grid points.}

\subsubsection{Spectral Monotonicity}
\RV{To compare different lag times on a common time horizon, fix $T>0$ and
set $\Delta_k=T/k$ for $k\geq1$. For $k=1$, the chain makes one transition
over $[0,T]$ and is exact by Proposition~\ref{prop:one_step}.
Increasing $k$ introduces more Markov transitions over the same interval and
therefore discards the unresolved memory more frequently. Define}
\begin{equation}\label{Hmc_T}
    \widetilde{\mathbf{H}}_k(T)
    :=
    \beta\mathbf{B}_0^{1/2}
    \widetilde{\mathbf{R}}_k(T)
    \mathbf{B}_0^{1/2}=\beta\mathbf{B}_0^{1/2}
    \widetilde{\mathbf{R}}(k\Delta_k;\Delta_k)
    \mathbf{B}_0^{1/2}=\mathbf{H}(\Delta_k)^{k}.
\end{equation}

\RV{For a symmetric matrix $\mathbf M\in\mathbb R^{n\times n}$, let
$\lambda_1(\mathbf M)\ge\cdots\ge\lambda_n(\mathbf M)$ denote its eigenvalues
in decreasing order.}

\begin{theorem}[Spectral monotonicity]\label{thm:monotonicity}
For $j=1,\ldots,n$,
\begin{equation}\label{monotonicity}
    \lambda_j\!\left(\widetilde{\mathbf{H}}_{k+1}(T)\right)
    \le
    \lambda_j\!\left(\widetilde{\mathbf{H}}_k(T)\right),
    \qquad k\ge1.
\end{equation}
Moreover, when $k \to \infty$,
\begin{equation}\label{Hmc_limit}
    \widetilde{\mathbf{H}}_k(T)
    \longrightarrow
    e^{-\mathbf{B}_0T}.
\end{equation}
\end{theorem}

\RV{Theorem~\ref{thm:monotonicity}, proved in Appendix~\ref{app:statistical},
gives the following upper and lower spectral bounds:}
\begin{equation}\label{eig_bounds}
    \RV{\lambda_j\!\left(\mathbf H(T)\right)
    =\lambda_j\!\left(\widetilde{\mathbf H}_1(T)\right)
    \ge \lambda_j\!\left(\widetilde{\mathbf H}_k(T)\right)
    \ge \lambda_j\!\left(e^{-\mathbf B_0T}\right),
    \qquad j=1,\ldots,n.}
\end{equation}
\RV{Thus one transition reproduces the exact autocovariance at $T$, whereas
successive refinement monotonically lowers every normalized autocovariance
mode and approaches a limit.}

\begin{corollary}
\label{cor:acf_bounds}
\RV{If $n=1$, then}
\begin{equation}\label{Rbounds}
    R(T)
    =
    \widetilde{R}_1(T)
    \ge
    \widetilde{R}_2(T)
    \ge
    \cdots
    \ge
    R(0)e^{-B_0T},
\end{equation}
\RV{and the inequalities between consecutive terms are strict if $A_0>B_0$.
In any dimension,
$\operatorname{Tr}\widetilde{\mathbf H}_k(T)$,
$\|\widetilde{\mathbf H}_k(T)\|_2$,
$\|\widetilde{\mathbf H}_k(T)\|_{\mathrm F}$, and
$\det\widetilde{\mathbf H}_k(T)$ are nonincreasing in $k$.}
\end{corollary}

\begin{proof}
\RV{All statements follow from the ordered eigenvalue monotonicity in
Theorem~\ref{thm:monotonicity}. When $n=1$,
$\widetilde R_k(T)=R(0)\widetilde H_k(T)$.}
\end{proof}

\section{Analysis for A 2D System}\label{sec:2d}
\RV{The general autocovariance analysis gives structural matrix conclusions,
but the dependence on the lag time remains implicit. We therefore
specialize to a two-dimensional linear system with one slow and one fast
variable. Motivated by
\cite{hudson2024dynamical}, we set}
\begin{equation}\label{model2d}
\begin{gathered}
    \mathbf{K}
    =\begin{pmatrix}1&0\\0&\lambda\end{pmatrix},
    \qquad
    \bPhi=\begin{pmatrix}\cos\theta&-\sin\theta\end{pmatrix},
    \qquad
    \bPsi=\begin{pmatrix}\sin\theta&\cos\theta\end{pmatrix}.
\end{gathered}
\end{equation}
\RV{We assume $\lambda>1$ and write $c=\cos\theta$, $s=\sin\theta$,
$p(t)=e^{-t}$, and $q(t)=e^{-\lambda t}$. Throughout this section,
$c\neq0$. If $c=0$, the resolved dynamics is already closed and no
coarse-graining is needed. The initial condition is
$(\xi_0,\zeta_0)=(\mu_{\xi,0},\mu_{\zeta,0})$.}

\subsection{Error of Autocovariance} \label{sec:error2d}

\RV{We first examine how increasing the time-scale separation affects the
autocovariance error when the lag time is fixed. The following theorem gives
the leading-order error as $\lambda\to\infty$ for each fixed step $k$.
The exact formulas and proofs of Theorems~\ref{thm:fixed_lag}
and~\ref{thm:small_lag} are collected in Appendix~\ref{app:slow_fast}.}

\begin{theorem}[Fixed Time-lag]\label{thm:fixed_lag}
\RV{Under the setting of this section, fix $k\in\mathbb{N}$ with $k\ge1$
and $\Delta>0$. As $\lambda\to\infty$,}
\begin{equation}
    \frac{|\widetilde{R}(k\Delta)-R(k\Delta)|}{R(k\Delta)}
    =
    \frac{|1-k|\tan^2\theta}{\lambda}
    + O(\lambda^{-2}).
    \label{error2d_fixed}
\end{equation}
\end{theorem}

\RV{We next allow the lag time to decrease as the time-scale separation
increases. The following theorem describes how small $\Delta$ can be while
still obtaining a vanishing autocovariance error at a fixed observation time.}

\begin{theorem}[\RV{Small-lag asymptotics}]\label{thm:small_lag}
\RV{Under the setting of this section, fix $\tau>0$. Let
$\Delta\to0$ and $\lambda\to\infty$, and choose
$k=k(\Delta)\in\mathbb{N}$ such that $k\Delta\to\tau$.
Depending on the relative scaling of $\lambda$ and $\Delta$,
the following asymptotic expansions hold:}
\begin{itemize}
    \item \RV{If $\lambda\Delta\to0$, then}
    \begin{equation}
        \RV{
        \frac{|\widetilde{R}(k\Delta)-R(k\Delta)|}{R(k\Delta)}
        =1-e^{-\tan^2\theta\,\tau}
        +O\!\left(|k\Delta-\tau|+\lambda\Delta+\lambda^{-1}\right).}
        \label{error2d_small}
    \end{equation}

    \item \RV{If $\lambda\Delta\to\ell\in(0,\infty)$, then}
    \begin{equation}
        \RV{
        \frac{|\widetilde{R}(k\Delta)-R(k\Delta)|}{R(k\Delta)}
        =1-e^{-\tan^2\theta\,\tau(1-e^{-\ell})/\ell}
        +O\!\left(|k\Delta-\tau|+|\lambda\Delta-\ell|+\Delta\right).}
        \label{error2d_mid}
    \end{equation}

    \item \RV{If $\lambda\Delta\to\infty$, then}
    \begin{equation}
        \RV{
        \frac{|\widetilde{R}(k\Delta)-R(k\Delta)|}{R(k\Delta)}
        =\frac{\tan^2\theta\,\tau}{\lambda\Delta}
        +O\!\left(
        \frac{|k\Delta-\tau|}{\lambda\Delta}
        +\lambda^{-1}+\frac{1}{(\lambda\Delta)^2}
        \right).}
        \label{error2d_limit}
    \end{equation}
\end{itemize}
\end{theorem}

\RV{The estimates \eqref{error2d_small}--\eqref{error2d_limit} show that the
lag time should not be taken too small. For fixed $\tan\theta\ne0$,
convergence requires $\Delta$ to remain large compared with the fast
relaxation time $\lambda^{-1}$. Taking $\Delta$ of order $\lambda^{-1}$ or
smaller leaves a nonvanishing error, even as the time-scale separation
tends to infinity.}

\subsection{Choice of Lag Time}\label{sec:lag}
Theorem~\ref{thm:finite_time}, proved in Appendix~\ref{app:slow_fast},
gives a finite-horizon error bound and a direct criterion for choosing $\Delta$.

\begin{theorem}[Finite-horizon autocovariance error]
\label{thm:finite_time}
Under the setting of this section, fix $T>0$. For every $0<\Delta\le T$,
\begin{equation}\label{error2d_bound}
    0
    \le
    \sup_{1\le k\le \lfloor T/\Delta\rfloor}
    \frac{R(k\Delta)-\widetilde{R}(k\Delta)}{R(k\Delta)}
    \le
    1-
    \exp\left[
    -\frac{T\tan^2\theta}{\lambda\Delta}
    \left(1-e^{-(\lambda-1)\Delta}\right)
    \right].
\end{equation}
\end{theorem}

\begin{corollary}[Choice of lag time]\label{cor:lag}
Fix $T>0$ and $0<\varepsilon<1$. For $0<\Delta\le T$, the bound
\begin{equation}\label{lag_choice}
    \Delta
    \ge
    \frac{T\tan^2\theta}
    {\lambda[-\log(1-\varepsilon)]}
\end{equation}
is sufficient for the relative autocovariance error in
\eqref{error2d_bound} to be at most $\varepsilon$ on
$\{k\Delta:k=1,\ldots,\lfloor T/\Delta\rfloor\}$.
\end{corollary}


Condition \eqref{lag_choice} is a sufficient lower bound on
the lag time. Figure~\ref{fig:Ex1_lag} illustrates the error and upper bound
in Theorem~\ref{thm:finite_time}. The upper range of $\Delta$ is determined separately by the desired
resolution of the slow dynamics and the prediction horizon.

\section{Comparison with Continuous-Time Markovian Models}\label{sec:comparison}
In this section, we compare the finite-lag Markov chain with existing
continuous-time Markovian effective models. The models considered below are
formulated as SDEs for the resolved variables. In the present linear Gaussian
setting, these SDEs can be written explicitly,
which allows a direct comparison of their transition laws and equilibrium
autocovariance functions. We focus on the following two representative
approaches.
\begin{itemize}
    \item \textbf{Approach~$\mathtt{A\text{-}1}$: projected or averaged Markovian dynamics.}
    A classical way to construct a Markovian effective dynamics is to design
    an infinitesimal SDE for the resolved variables by eliminating the
    unresolved variables at the level of the drift, diffusion, or generator.
    In the linear setting considered here, many such constructions
    lead to the same effective Ornstein--Uhlenbeck dynamics. We refer to this
    common effective model as Approach~$\mathtt{A\text{-}1}$:
    \begin{equation*}
      \diff\widetilde{\boldsymbol{\xi}}^{\mathtt{A\text{-}1}}_t
      =
      -\mathbf{B}_0\widetilde{\boldsymbol{\xi}}^{\mathtt{A\text{-}1}}_t \diff t
      +
      \sqrt{2\beta^{-1}}\,\diff\mathbf{w}_t .
    \end{equation*}
    \RV{More concretely, conditional-expectation dynamics derives this model by
    averaging the full coefficients over the equilibrium conditional
    distribution on each level set of the coarse-graining map
    \cite{legoll2010effective}.}
    \item \textbf{Approach~$\mathtt{A\text{-}2}$: memory-corrected Markovian dynamics.}
    Hudson and Li proposed a corrected Markovian dynamics for coarse-grained
    linear SDEs \cite{hudson2024dynamical}. This model modifies the projected
    dynamics by incorporating a Markovian correction motivated by the memory
    effect of unresolved variables. We refer to this model as Approach~$\mathtt{A\text{-}2}$:
    \begin{equation*}
      \diff\widetilde{\boldsymbol{\xi}}^{\mathtt{A\text{-}2}}_t
      =
      -\mathbf{C}\mathbf{B}_0
      \widetilde{\boldsymbol{\xi}}^{\mathtt{A\text{-}2}}_t \diff t
      +
      \sqrt{2\beta^{-1}}\,\mathbf{C}^{1/2}\diff\mathbf{w}_t,
      \qquad
      \mathbf{C}
      =
      \left(
      \mathbf{I}_n+\boldsymbol{\alpha}\mathbf{A}_1^{-2}\boldsymbol{\alpha}^{\top}
      \right)^{-1}.
    \end{equation*}
\end{itemize}

\subsection{A Unified Representation}\label{sec:unified}
We now place the finite-lag chain and the above two continuous-time Markovian approximations in the same framework. All three models preserve the same equilibrium distribution \cite{hudson2024dynamical}, but they represent temporal correlations in different ways. Normalizing the equilibrium covariance separates these two
properties and reduces the comparison to the choice of a single relaxation matrix.

For a resolved variable with equilibrium distribution
$\mathcal{N}(\mathbf{0},\beta^{-1}\mathbf{B}_0^{-1})$, we introduce the normalized resolved coordinate and its equilibrium autocovariance:
\begin{equation}\label{normalized}
    \mathbf{y}_t
    :=
    \beta^{1/2}\mathbf{B}_0^{1/2}\boldsymbol{\xi}_t,
    \qquad
    \mathbf{H}(t)
    :=
    \RV{\mathbb{E}_{\pi_{\mathbf x}}[\mathbf{y}_t\mathbf{y}_0^\top]}
    =
    \beta\mathbf{B}_0^{1/2}\mathbf{R}(t)\mathbf{B}_0^{1/2}.
\end{equation}
Thus $\mathbf{y}_t$ has equilibrium distribution
$\mathcal{N}(\mathbf{0},\mathbf{I}_n)$, and $\mathbf{H}(t)$ is the normalized
equilibrium autocovariance of the exact resolved dynamics: it measures the
linear dependence between the resolved states at equilibrium separated by a
lag $t$, after removing the scale set by their equilibrium covariance.

We consider an Ornstein--Uhlenbeck (OU) process preserving
$\mathcal{N}(\mathbf{0},\mathbf{I}_n)$ of the form
\begin{equation}\label{OU}
    \diff\widetilde{\mathbf{y}}_t
    =
    -\mathbf{L}\widetilde{\mathbf{y}}_t\,\diff t
    +
    \sqrt{2\mathbf{L}}\,\diff\mathbf{w}_t,
    \qquad
    \mathbf{L}=\mathbf{L}^\top\succ\mathbf{0}.
\end{equation}
Here $\mathbf{L}\in\mathbb{R}^{n\times n}$ is a constant symmetric positive
definite matrix, and $\sqrt{2\mathbf{L}}$ denotes its symmetric positive definite square root. At equilibrium, the corresponding normalized equilibrium autocovariance is given by
\begin{equation}\label{Hou}
    \widetilde{\mathbf{H}}(t;\mathbf{L})
    :=
    \RV{\mathbb{E}_{\widetilde{\mathbf{y}}_0\sim
    \mathcal{N}(\mathbf{0},\mathbf{I}_n)}}
    \left[
    \widetilde{\mathbf{y}}_t\widetilde{\mathbf{y}}_0^\top
    \right]
    =
    e^{-\mathbf{L}t}
    \mathbb{E}
    \left[
    \widetilde{\mathbf{y}}_0\widetilde{\mathbf{y}}_0^\top
    \right]
    =
    e^{-\mathbf{L}t}.
\end{equation}

Unless $\mathbf{H}(t)$ has the semigroup property, no constant matrix
$\mathbf{L}$ can make $\widetilde{\mathbf{H}}(t;\mathbf{L})=\mathbf{H}(t)$ for every $t\ge0$ since $\widetilde{\mathbf{H}}(t)$ naturally has semigroup property. In this section, we illustrate that the below models can be unified by selecting $\mathbf{L}$ preserving different features of the exact correlation curve:
\begin{itemize}[itemsep=0.15em,topsep=0.25em]
    \item \textbf{A Continuous Embedding of the Finite-Lag Chain}: Finite-lag matching preserves the exact normalized autocovariance at the lag time $\Delta$: $$\widetilde{\mathbf{H}}^{\mathtt{MC}}(\Delta)=\mathbf{H}(\Delta)\text{, with fixed }\Delta.$$
    The resulting SDE is a continuous embedding of the finite-lag Markov chain.
    \item \textbf{Approach~$\mathtt{A\text{-}1}$}: Initial-slope matching preserves the infinitesimal relaxation rate: $$\left(\widetilde{\mathbf{H}}^{\mathtt{A\text{-}1}}\right)'(0)=\mathbf{H}'(0).$$
    \item \textbf{Approach~$\mathtt{A\text{-}2}$}: Time-integral matching preserves the total correlation carried over all lags:
    $$\int_0^\infty\widetilde{\mathbf{H}}^{\mathtt{A\text{-}2}}(t)\,\diff t=\int_0^\infty\mathbf{H}(t)\,\diff t.$$
\end{itemize}

\subsection{A Continuous Embedding of the Finite-Lag Chain}\label{sec:embedding}
For fixed $\Delta>0$, define
\begin{equation}\label{Lmc}
    \mathbf{L}^{\mathtt{MC}}_{\Delta}
    :=
    -\frac{1}{\Delta}\log\mathbf{H}(\Delta).
\end{equation}
Since $\mathbf{0}\prec\mathbf{H}(\Delta)\prec\mathbf{I}_n$, the principal
matrix logarithm in \eqref{Lmc} is well defined
and $\mathbf{L}^{\mathtt{MC}}_{\Delta}\succ\mathbf{0}$. This choice is
characterized by the finite-lag matching condition
\begin{equation}\label{Hmatch}
    \widetilde{\mathbf{H}}^{\mathtt{MC}}(t)
    :=
    e^{-\mathbf{L}^{\mathtt{MC}}_{\Delta}t},
    \qquad
    \widetilde{\mathbf{H}}^{\mathtt{MC}}(\Delta)
    =
    \mathbf{H}(\Delta).
\end{equation}
Thus the embedded dynamics reproduces the exact equilibrium correlation at
the lag time $\Delta$. This is the relevant finite-time statistic for a model
whose states are advanced on the grid $\{k\Delta\}_{k\geq0}$.
At the skeleton times, this notation agrees with Section~\ref{sec:acf}:
\begin{equation}\label{Hmc_grid}
    \widetilde{\mathbf{H}}^{\mathtt{MC}}(k\Delta)
    =
    \widetilde{\mathbf{H}}(k\Delta)
    =
    \mathbf{H}(\Delta)^k.
\end{equation}
Note that in general, $\mathbf{H}(\Delta)^k\neq\mathbf{H}(k\Delta)$, which reflects the failure of the normalized autocovariance of the resolved process to satisfy the semigroup property.

\RV{For random elements $X$ and $Y$, we write
$X\overset{\mathrm{law}}{=}Y$ if they have the same probability law.}

\begin{proposition}[Continuous embedding on Normalized Coordinates]\label{prop:embedding}
\RV{Fix $\Delta>0$. Let $\{\widetilde{\mathbf{y}}_t\}_{t\geq0}$ solve
\eqref{OU} with
$\mathbf{L}=\mathbf{L}^{\mathtt{MC}}_{\Delta}$ and
$\widetilde{\mathbf{y}}_0\overset{\mathrm{law}}{=}
\beta^{1/2}\mathbf{B}_0^{1/2}\widetilde{\boldsymbol{\xi}}_0$. Then}
\[
    \RV{
    \bigl(\widetilde{\mathbf{y}}_{k\Delta}\bigr)_{k\geq0}
    \overset{\mathrm{law}}{=}
    \bigl(\beta^{1/2}\mathbf{B}_0^{1/2}
    \widetilde{\boldsymbol{\xi}}_{k\Delta}\bigr)_{k\geq0}.}
\]
\end{proposition}

\begin{proof}
\RV{Set $\mathbf{L}=\mathbf{L}^{\mathtt{MC}}_{\Delta}$. The solution of
\eqref{OU} over one interval of length $\Delta$ satisfies}
\begin{equation}\label{OU_step}
    \RV{
    \widetilde{\mathbf{y}}_{t+\Delta}
    =e^{-\mathbf{L}\Delta}\widetilde{\mathbf{y}}_t
    +\int_t^{t+\Delta}
    e^{-\mathbf{L}(t+\Delta-s)}\sqrt{2\mathbf{L}}\,\diff\mathbf{w}_s.}
\end{equation}
\RV{Since the stochastic integral is centered Gaussian and independent of
$\widetilde{\mathbf{y}}_t$, for every $\mathbf{y}\in\mathbb{R}^n$,}
\begin{align}\label{OU_kernel}
    \RV{
    \widetilde{\mathbf{y}}_{t+\Delta}
    \mid \widetilde{\mathbf{y}}_t=\mathbf{y}
    }
    &\RV{\sim\mathcal{N}\!\left(
    e^{-\mathbf{L}\Delta}\mathbf{y},
    \int_0^\Delta e^{-\mathbf{L}s}(2\mathbf{L})e^{-\mathbf{L}s}\,\diff s
    \right)} \RV{=}\RV{\mathcal{N}\!\left(
    \mathbf{H}(\Delta)\mathbf{y},
    \mathbf{I}_n-\mathbf{H}(\Delta)^2
    \right),}
\end{align}
\RV{where $e^{-\mathbf{L}\Delta}=\mathbf{H}(\Delta)$ and
$\int_0^\Delta 2\mathbf{L}e^{-2\mathbf{L}s}\,\diff s
=\mathbf{I}_n-e^{-2\mathbf{L}\Delta}$. Applying the same normalization to
the Markov-chain recursion \eqref{MC} gives}
\begin{align}\label{MC_normalized}
    \RV{
    \beta^{1/2}\mathbf{B}_0^{1/2}
    \widetilde{\boldsymbol{\xi}}_{(k+1)\Delta}
    \mid
    \beta^{1/2}\mathbf{B}_0^{1/2}
    \widetilde{\boldsymbol{\xi}}_{k\Delta}=\mathbf{y}
    }
    &\RV{\sim\mathcal{N}\!\left(
    \mathbf{B}_0^{1/2}\mathbf{F}_\infty\mathbf{B}_0^{-1/2}\mathbf{y},
    \beta\mathbf{B}_0^{1/2}\mathbf{G}_\infty\mathbf{B}_0^{1/2}
    \right).}
\end{align}
\RV{By \eqref{H} and \eqref{Gidentity},}
\begin{align}\label{MC_parameters}
    \RV{\mathbf{B}_0^{1/2}\mathbf{F}_\infty\mathbf{B}_0^{-1/2}}
    &=\RV{\mathbf{H}(\Delta),} \notag\\
    \RV{\beta\mathbf{B}_0^{1/2}\mathbf{G}_\infty\mathbf{B}_0^{1/2}}
    &=\RV{\mathbf{I}_n-\mathbf{H}(\Delta)\mathbf{H}(\Delta)^\top
    =\mathbf{I}_n-\mathbf{H}(\Delta)^2.}
\end{align}
\RV{This completes the proof.}
\end{proof}

Returning to the original resolved coordinate, define
\begin{equation}\label{Dmc}
    \mathbf{D}^{\mathtt{MC}}_{\Delta}
    :=
    \mathbf{B}_0^{-1/2}
    \mathbf{L}^{\mathtt{MC}}_{\Delta}
    \mathbf{B}_0^{1/2}
    =
    -\frac{1}{\Delta}\log\mathbf{F}_\infty.
\end{equation}
A continuous embedding of the finite-lag chain is therefore given by
\begin{equation}\label{embedding}
    \diff\widetilde{\boldsymbol{\xi}}^{\mathtt{MC}}_t
    =
    -\mathbf{D}^{\mathtt{MC}}_{\Delta}
    \widetilde{\boldsymbol{\xi}}^{\mathtt{MC}}_t\,\diff t
    +
    \sqrt{2\beta^{-1}}\,
    \mathbf{B}_0^{-1/2}
    \left(\mathbf{L}^{\mathtt{MC}}_{\Delta}\right)^{1/2}
    \,\diff\mathbf{w}_t.
\end{equation}

\subsection{Comparison with Approach 1}
A classical Markovian approximation is obtained by eliminating the
unresolved variables at the level of the drift, diffusion, or generator. In
\RV{the present linear setting, conditional-expectation dynamics
\cite{legoll2010effective}, projected-generator and Dirichlet-form methods
\cite{nuske2021spectral} lead to the same effective
Ornstein--Uhlenbeck dynamics:}
\begin{equation*}\tag{A-1}\label{effA1}
    \diff\widetilde{\boldsymbol{\xi}}^{\mathtt{A\text{-}1}}_t
    =
    -\mathbf{B}_0
    \widetilde{\boldsymbol{\xi}}^{\mathtt{A\text{-}1}}_t\,\diff t
    +
    \sqrt{2\beta^{-1}}\,\diff\mathbf{w}_t.
\end{equation*}
Writing
$\widetilde{\mathbf{y}}^{\mathtt{A\text{-}1}}_t
:=\beta^{1/2}\mathbf{B}_0^{1/2}
\widetilde{\boldsymbol{\xi}}^{\mathtt{A\text{-}1}}_t$,
Approach~$\mathtt{A\text{-}1}$ takes the normalized form
\begin{equation}\label{A1_normalized}
    \diff\widetilde{\mathbf{y}}^{\mathtt{A\text{-}1}}_t
    =
    -\mathbf{B}_0
    \widetilde{\mathbf{y}}^{\mathtt{A\text{-}1}}_t\,\diff t
    +
    \sqrt{2\mathbf{B}_0}\,\diff\mathbf{w}_t.
\end{equation}
Thus its normalized drift in \eqref{OU} is
\begin{equation}\label{La1}
    \mathbf{L}^{\mathtt{A\text{-}1}}
    =
    \mathbf{B}_0.
\end{equation}
\RV{Its normalized equilibrium autocovariance is
$\widetilde{\mathbf{H}}^{\mathtt{A\text{-}1}}(t)=e^{-\mathbf{B}_0t}$.}
Since $\mathbf{H}'(0)=-\mathbf{B}_0$ by
\eqref{U}, this choice satisfies the local matching
condition
\begin{equation}\label{slope_match}
    \left(
    \widetilde{\mathbf{H}}^{\mathtt{A\text{-}1}}
    \right)'(0)
    =
    -\mathbf{L}^{\mathtt{A\text{-}1}}
    =
    -\mathbf{B}_0
    =
    \mathbf{H}'(0).
\end{equation}
Approach~$\mathtt{A\text{-}1}$ is constructed from an infinitesimal drift or
generator, so matching the instantaneous relaxation rate is the corresponding
local criterion. It preserves the first-order behavior of the exact
autocovariance at short lags, but does not impose a finite-lag or integrated
matching condition.

\begin{lemma}[Local correspondence]\label{lem:local}
As $\Delta\to0$,
\begin{equation}\label{Lmc_limit}
    \mathbf{L}^{\mathtt{MC}}_{\Delta}
    \longrightarrow
    \mathbf{L}^{\mathtt{A\text{-}1}}
    =
    \mathbf{B}_0.
\end{equation}
\end{lemma}

The local correspondence of the normalized drifts in Lemma~\ref{lem:local},
proved in Appendix~\ref{app:local}, also gives a global process-level limit.

\begin{theorem}[Small-lag limit toward Approach 1]
\label{thm:A1_limit}
Let $\widetilde{\mathbf{y}}^{\mathtt{MC},\Delta}$ be the continuous embedding
of the Markov chain. Couple it with
$\widetilde{\mathbf{y}}^{\mathtt{A\text{-}1}}$ using the same Brownian motion
and the same initial condition, with finite second moment. Then, for every
$T>0$,
\begin{equation}\label{path_limit}
    \mathbb{E}\left[
    \sup_{0\le t\le T}
    \left\|
    \widetilde{\mathbf{y}}^{\mathtt{MC},\Delta}_t
    -\widetilde{\mathbf{y}}^{\mathtt{A\text{-}1}}_t
    \right\|^2
    \right]
    \longrightarrow0
    \qquad\text{as }\Delta\to0.
\end{equation}
\end{theorem}

\begin{proof}
Lemma~\ref{lem:local} and the continuity of the principal
matrix square root give
\begin{equation}\label{coeff_limit}
    \mathbf{L}^{\mathtt{MC}}_{\Delta}\longrightarrow\mathbf{B}_0,
    \qquad
    \sqrt{2\mathbf{L}^{\mathtt{MC}}_{\Delta}}
    \longrightarrow\sqrt{2\mathbf{B}_0}.
\end{equation}
The drift and diffusion coefficients are globally Lipschitz, with uniformly
bounded Lipschitz and growth constants for sufficiently small $\Delta$.
The standard continuous-dependence estimate for SDE solutions therefore
gives \eqref{path_limit}.
\RV{See \cite[Theorem~4.5.6, the subsequent remark, and
Section~4.8]{kloeden1992numerical}.}
\end{proof}

\RV{For a fixed full system, Theorem~\ref{thm:A1_limit} also implies
the fixed-time autocovariance limit \eqref{Hmc_limit}.}

\subsection{Comparison with Approach 2}
Hudson and Li proposed the memory-corrected Markovian dynamics
\cite{hudson2024dynamical}
\begin{equation*}\tag{A-2}\label{effA2}
    \diff\widetilde{\boldsymbol{\xi}}^{\mathtt{A\text{-}2}}_t
    =
    -\mathbf{C}\mathbf{B}_0
    \widetilde{\boldsymbol{\xi}}^{\mathtt{A\text{-}2}}_t\,\diff t
    +
    \sqrt{2\beta^{-1}}\,\mathbf{C}^{1/2}\diff\mathbf{w}_t,
    \qquad
    \mathbf{C}
    =
    \left(
    \mathbf{I}_n+
    \boldsymbol{\alpha}\mathbf{A}_1^{-2}\boldsymbol{\alpha}^{\top}
    \right)^{-1}.
\end{equation*}
Writing
$\widetilde{\mathbf{y}}^{\mathtt{A\text{-}2}}_t
:=\beta^{1/2}\mathbf{B}_0^{1/2}
\widetilde{\boldsymbol{\xi}}^{\mathtt{A\text{-}2}}_t$,
its normalized drift is
\begin{equation}\label{La2}
    \mathbf{L}^{\mathtt{A\text{-}2}}
    :=
    \mathbf{B}_0^{1/2}\mathbf{C}\mathbf{B}_0^{1/2}.
\end{equation}
The normalized diffusion in \eqref{effA2} has covariance
$2\mathbf{L}^{\mathtt{A\text{-}2}}$, so the normalized process is
statistically equivalent to \eqref{OU} with
$\mathbf{L}=\mathbf{L}^{\mathtt{A\text{-}2}}$.
\RV{Its normalized equilibrium autocovariance is therefore
$\widetilde{\mathbf{H}}^{\mathtt{A\text{-}2}}(t)
=e^{-\mathbf{L}^{\mathtt{A\text{-}2}}t}$.}
As proved below, this choice is characterized by
\begin{equation}\label{integral_match}
    \int_0^\infty
    \widetilde{\mathbf{H}}^{\mathtt{A\text{-}2}}(t)\,\diff t
    =
    \left(\mathbf{L}^{\mathtt{A\text{-}2}}\right)^{-1}
    =
    \int_0^\infty\mathbf{H}(t)\,\diff t.
\end{equation}
Thus Approach~$\mathtt{A\text{-}2}$ matches the total equilibrium correlation
accumulated over all lags. This integral gives the matching condition a direct long-time
interpretation.

\begin{proposition}[Integrated autocovariance matching]
\label{prop:integral}
Approach~$\mathtt{A\text{-}2}$ satisfies
\begin{equation}\label{Rintegral_match}
    \int_0^\infty
    \widetilde{\mathbf{R}}^{\mathtt{A\text{-}2}}(t)\,\diff t
    =
    \int_0^\infty\mathbf{R}(t)\,\diff t.
\end{equation}
Equivalently,
\eqref{integral_match} holds.
\end{proposition}

\begin{proof}
The exact autocovariance gives
\begin{equation}\label{Rintegral}
    \int_0^\infty\mathbf{R}(t)\,\diff t
    =
    \beta^{-1}[\mathbf{A}^{-2}]_{\xi\xi}.
\end{equation}
The $\xi$-$\xi$ block of $\mathbf{A}^{-1}$ is $\mathbf{B}_0^{-1}$, while its
$\xi$-$\zeta$ block is
$-\mathbf{B}_0^{-1}\boldsymbol{\alpha}\mathbf{A}_1^{-1}$. Hence
\begin{align}\label{Ainv2}
    [\mathbf{A}^{-2}]_{\xi\xi}
    &=
    \mathbf{B}_0^{-2}
    +
    \mathbf{B}_0^{-1}
    \boldsymbol{\alpha}\mathbf{A}_1^{-2}\boldsymbol{\alpha}^{\top}
    \mathbf{B}_0^{-1}\notag\\
    &=
    \mathbf{B}_0^{-1}\mathbf{C}^{-1}\mathbf{B}_0^{-1}.
\end{align}
Approach~$\mathtt{A\text{-}2}$ has equilibrium autocovariance
\begin{equation}\label{Ra2}
    \widetilde{\mathbf{R}}^{\mathtt{A\text{-}2}}(t)
    =
    \beta^{-1}\mathbf{B}_0^{-1}e^{-\mathbf{B}_0\mathbf{C}t},
\end{equation}
and hence
\begin{equation}\label{Ra2_integral}
    \int_0^\infty
    \widetilde{\mathbf{R}}^{\mathtt{A\text{-}2}}(t)\,\diff t
    =
    \beta^{-1}\mathbf{B}_0^{-1}
    (\mathbf{B}_0\mathbf{C})^{-1}
    =
    \beta^{-1}\mathbf{B}_0^{-1}\mathbf{C}^{-1}\mathbf{B}_0^{-1}.
\end{equation}
Conjugating
\eqref{Rintegral_match} by
$\beta^{1/2}\mathbf{B}_0^{1/2}$ gives
\eqref{integral_match}.
\end{proof}

\RV{For a fixed $\Delta>0$ and identical initial laws, equality in law of
the two sampled processes requires
$\mathbf{L}^{\mathtt{MC}}_{\Delta}=\mathbf{L}^{\mathtt{A\text{-}2}}$.}
Equivalently,
\begin{equation}\label{A2_match}
    \mathbf{H}(\Delta)
    =
    e^{-\mathbf{L}^{\mathtt{A\text{-}2}}\Delta},
    \qquad\text{or}\qquad
    \mathbf{F}_\infty
    =
    e^{-\mathbf{C}\mathbf{B}_0\Delta}.
\end{equation}
\RV{For $n>1$, a scalar $\Delta$ need not satisfy
\eqref{A2_match}.} When the resolved variable is scalar,
an exact matching lag always exists.

\RV{Suppose $n=1$. Since $\mathbf{A}$ is symmetric positive definite, its
spectral decomposition and \eqref{U} give}
\begin{equation}\label{Hscalar}
    \RV{H(t)=\sum_{j=1}^Jp_je^{-\lambda_jt},}
\end{equation}
\RV{where $\lambda_1,\ldots,\lambda_J$ are the distinct decay rates appearing
in $H$, and $p_j>0$ with $\sum_{j=1}^Jp_j=1$.}

\begin{proposition}[Exact matching with a scalar resolved variable]
\label{prop:scalar_match}
\RV{Suppose $n=1$, with $J$ defined in \eqref{Hscalar}. If $J\ge2$, there exists a unique $\Delta^*>0$ such that}
\begin{equation}\label{Lmatch}
    L^{\mathtt{MC}}_{\Delta^*}
    =
    L^{\mathtt{A\text{-}2}}.
\end{equation}
Consequently,
\begin{equation}\label{law_match}
    \left\{
    \widetilde{\xi}^{\mathtt{MC}}_{k\Delta^*}
    \right\}_{k\ge0}
    \overset{\mathrm{law}}{=}
    \left\{
    \widetilde{\xi}^{\mathtt{A\text{-}2}}_{k\Delta^*}
    \right\}_{k\ge0}
\end{equation}
when initialized identically.
\RV{If $J=1$, \eqref{law_match} holds with $\Delta^*$ replaced by any
$\Delta>0$.}
\end{proposition}

\begin{proof}
\RV{By \eqref{Lmc}, \eqref{Hscalar} and \eqref{integral_match}}
\begin{equation}\label{scalar_rates}
    \RV{\begin{aligned}
    L^{\mathtt{MC}}_{\Delta}
    &=
    -\frac{1}{\Delta}
    \log\left(\sum_{j=1}^Jp_je^{-\lambda_j\Delta}\right),\\
    L^{\mathtt{A\text{-}2}}
    &=
    \left(\int_0^\infty H(t)\,\diff t\right)^{-1}
    =
    \left(\sum_{j=1}^Jp_j\int_0^\infty e^{-\lambda_jt}\,\diff t\right)^{-1}
    =
    \left(\sum_{j=1}^J\frac{p_j}{\lambda_j}\right)^{-1}.
    \end{aligned}}
\end{equation}
\RV{If $J=1$, then \eqref{Hscalar} gives $H(t)=e^{-\lambda_1t}$,
and \eqref{scalar_rates} yields
$L^{\mathtt{MC}}_{\Delta}=L^{\mathtt{A\text{-}2}}=\lambda_1$ for every
$\Delta>0$. Equation~\eqref{law_match} therefore holds for every $\Delta>0$
by Proposition~\ref{prop:embedding}.}

For $J\ge2$, the function
$-\log(\sum_jp_je^{-\lambda_j\Delta})$ is strictly concave and vanishes at
$\Delta=0$. Therefore $L^{\mathtt{MC}}_{\Delta}$ is continuous and strictly
decreasing, with
\begin{equation}\label{rate_limits}
    \RV{\lim_{\Delta\to0}L^{\mathtt{MC}}_{\Delta}
    =
    \sum_{j=1}^Jp_j\lambda_j,
    \qquad
    \lim_{\Delta\to\infty}L^{\mathtt{MC}}_{\Delta}
    =
    \min_{1\le j\le J}\lambda_j.}
\end{equation}
The value $L^{\mathtt{A\text{-}2}}$ is the weighted harmonic mean of the
active eigenvalues and lies strictly between the two limits in
\eqref{rate_limits}. The intermediate value theorem gives a
unique $\Delta^*$ satisfying \eqref{Lmatch}. Equation
\eqref{law_match} follows from
Proposition~\ref{prop:embedding}.
\end{proof}

\section{Numerical Examples}\label{sec:examples}
In this section, we present numerical examples to examine equilibrium convergence, autocovariance errors, and the comparison with existing continuous-time Markovian effective models. All statistics are evaluated from the moment relations and autocovariance formulas derived above, with $\beta=1$.

\subsection{A 2D Example}\label{ex:2d}
We test the two-dimensional model \eqref{model2d} from Section~\ref{sec:2d}, with $\theta=0.3$ unless noted.

\subsubsection{Convergence to Equilibrium}
\begin{figure}[!htbp]
    \centering
    \includegraphics[width=0.99\textwidth]{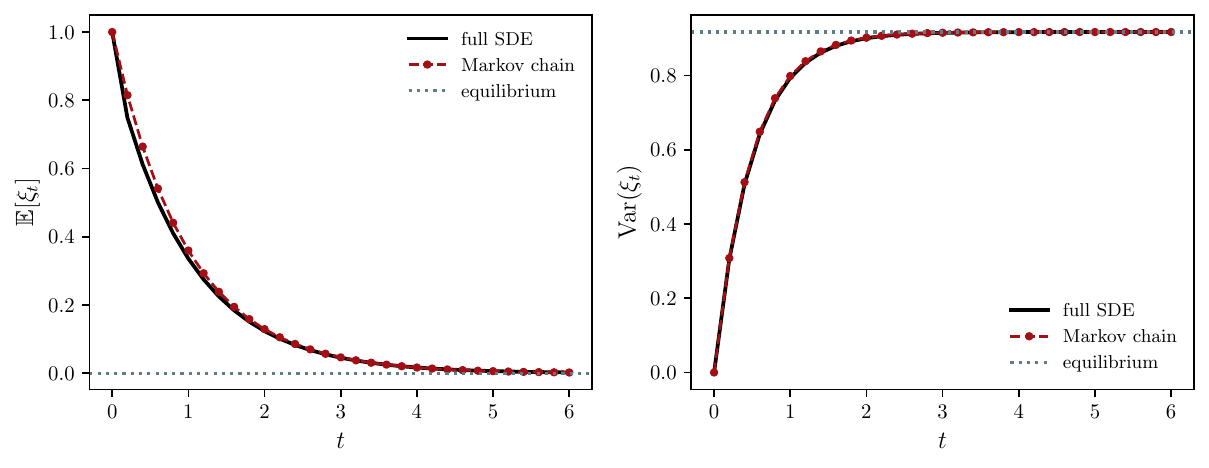}
    \caption{Resolved mean (left) and variance (right) for the full SDE and Markov chain, with $\lambda=20$, $\theta=0.3$, and $\Delta=0.2$.}
    \label{fig:Ex1_stats}
\end{figure}

We fix $\lambda=20$ and $\Delta=0.2$, and start from the deterministic
states $\mathbf{x}_0=(1,0)^\top$ and $\widetilde\xi_0=1$. Figure~\ref{fig:Ex1_stats} compares the resolved mean and variance for the effective model and the full SDE
up to $T=6$, computed using \eqref{moments}, \eqref{transMC_k},
and \eqref{covMC_limit}.
Both models approach the same equilibrium moments, despite differences in
their transient means. This agrees with Theorem~\ref{thm:equilibrium}.

\subsubsection{Relative Autocovariance Error}
\begin{figure}[!htbp]
    \centering
    \includegraphics[width=0.99\textwidth]{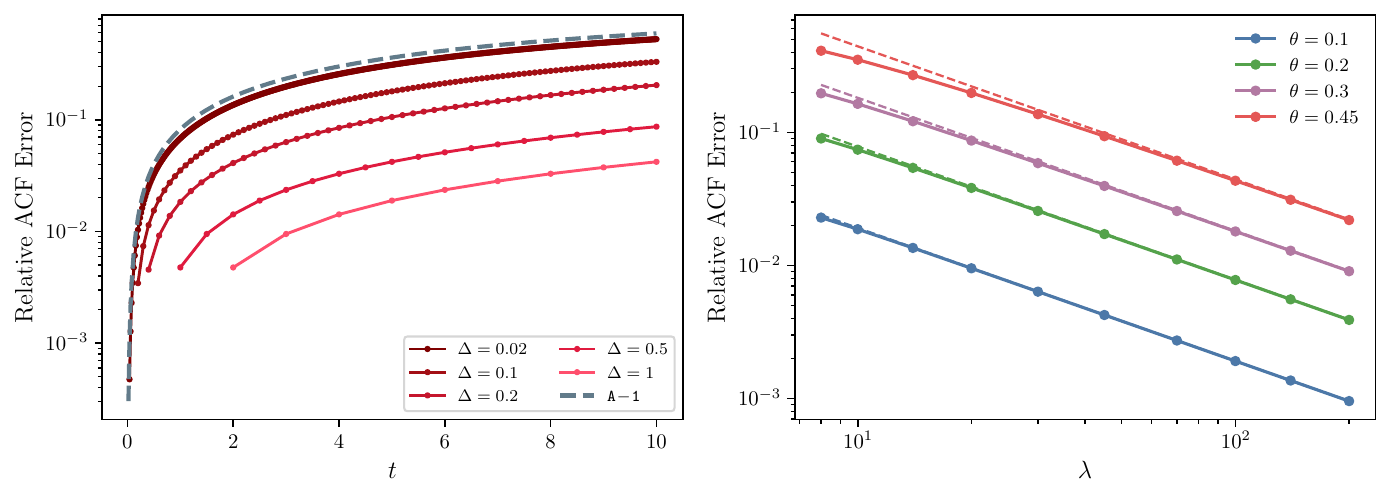}
    \caption{Relative ACF error in the 2D example. Left: varying $\Delta$ and steps $k$ at $\lambda=20$. Right: varying $\lambda$ at $\Delta=0.5$ and $k=20$.}
    \label{fig:Ex1_acf}
\end{figure}

The equilibrium autocovariances are evaluated using \eqref{exactR} and
\eqref{Rmc2d}. In Figure~\ref{fig:Ex1_acf}, the left panel fixes $\lambda=20$
and compares $\Delta=0.02,0.1,0.2,0.5,1$ up to $T=10$, with
Approach~$\mathtt{A\text{-}1}$ as a reference.
The errors approach the Approach~$\mathtt{A\text{-}1}$ curve as the lag decreases,
consistent with the limit discussed in Section~\ref{sec:comparison}.
The right panel fixes the lag time $\Delta=0.5$ and the number of transitions $k=20$, and varies
$\lambda\in[8,200]$ for $\theta=0.1,0.2,0.3,0.45$.
The errors approach the leading-order predictions of
Theorem~\ref{thm:fixed_lag}, confirming the inverse-$\lambda$ scaling.

\subsubsection{Lag-Time Selection}
\begin{figure}[!htbp]
    \centering
    \includegraphics[width=\textwidth]{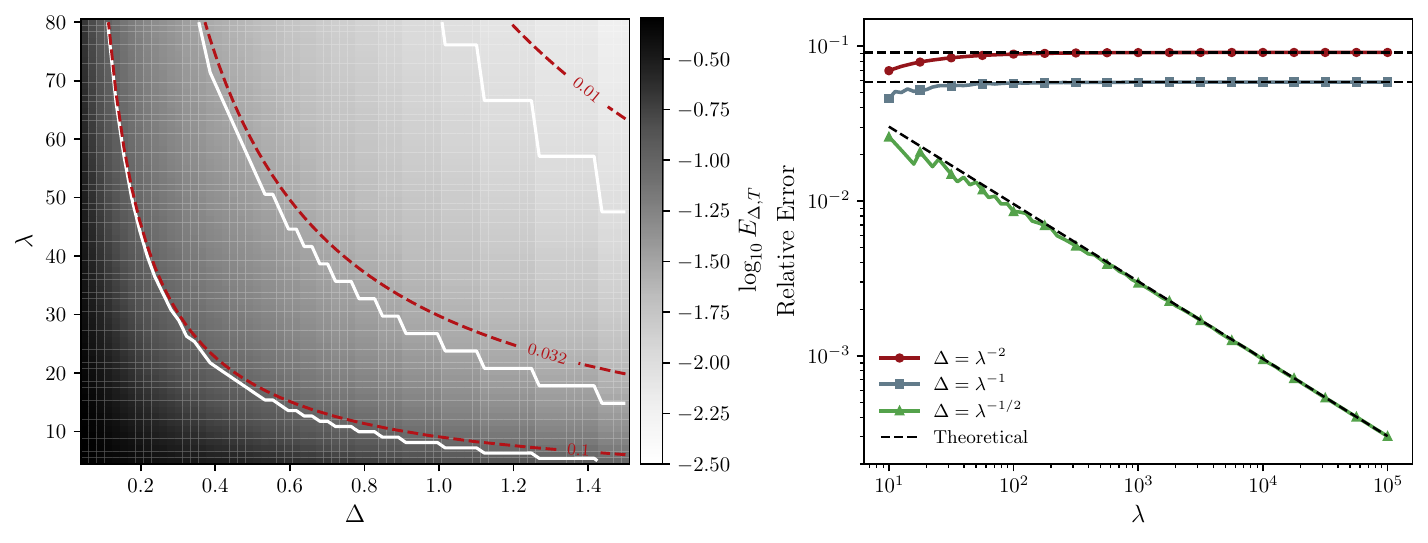}
    \caption{Relative ACF errors with $\theta=0.3$. Left: the finite-horizon error $E_{\Delta,T}$ and its upper bound from Theorem~\ref{thm:finite_time} for $T=10$, with contours shown in white and red, respectively. Right: relative error at $k=\lceil T/\Delta\rceil$ steps with $T=1$ for three small-lag scalings.}
    \label{fig:Ex1_lag}
\end{figure}

The left panel of Figure~\ref{fig:Ex1_lag} shows the finite-horizon error
$E_{\Delta,T}$ and its upper bound from Theorem~\ref{thm:finite_time},
with $T=10$, $\Delta\in[0.05,1.5]$, and $\lambda\in[5,80]$.
The computed errors remain below the theoretical upper bound throughout the sampled range.
The contours of the computed error (white) and the bound (red) exhibit
similar dependence on the lag and scale separation.

To examine the small-lag regimes, we fix $T=1$ and vary $\lambda\in[10,10^5]$, with
$\Delta=\lambda^{-2},\lambda^{-1},\lambda^{-1/2}$.
We evaluate the relative errors at $k=\lceil T/\Delta\rceil$ steps.
The right panel of Figure~\ref{fig:Ex1_lag} compares these relative errors from
\eqref{Rratio} with the predictions of Theorem~\ref{thm:small_lag}.
The first two scalings give nonzero limiting errors, as in
\eqref{error2d_small} and \eqref{error2d_mid}.
For the third scaling, the error decays linearly on log--log axes,
in agreement with \eqref{error2d_limit}.
Thus the three regimes distinguish convergence of the autocovariance
from a persistent relative error.

\subsection{A 40D-to-3D Example}\label{ex:chain}
We consider a harmonic chain with drift matrix
$\mathbf A_\ell=\kappa\mathbf I_{40}+\ell\mathbf L_{40}$, where $\kappa,\ell>0$ and
\[
\mathbf L_{40}=\begin{pmatrix}
1 & -1 & & & \\
-1 & 2 & -1 & & \\
& \ddots & \ddots & \ddots & \\
& & -1 & 2 & -1 \\
& & & -1 & 1
\end{pmatrix}
\]
is the path-graph Laplacian.
We take $q_{10}$, $q_{20}$, and $q_{30}$ as the resolved variables. We evaluate the autocovariances using \eqref{R} and \eqref{Rmc},
and measure the relative error in the matrix $2$-norm.

\begin{figure}[!htbp]
    \centering
    \includegraphics[width=\textwidth]{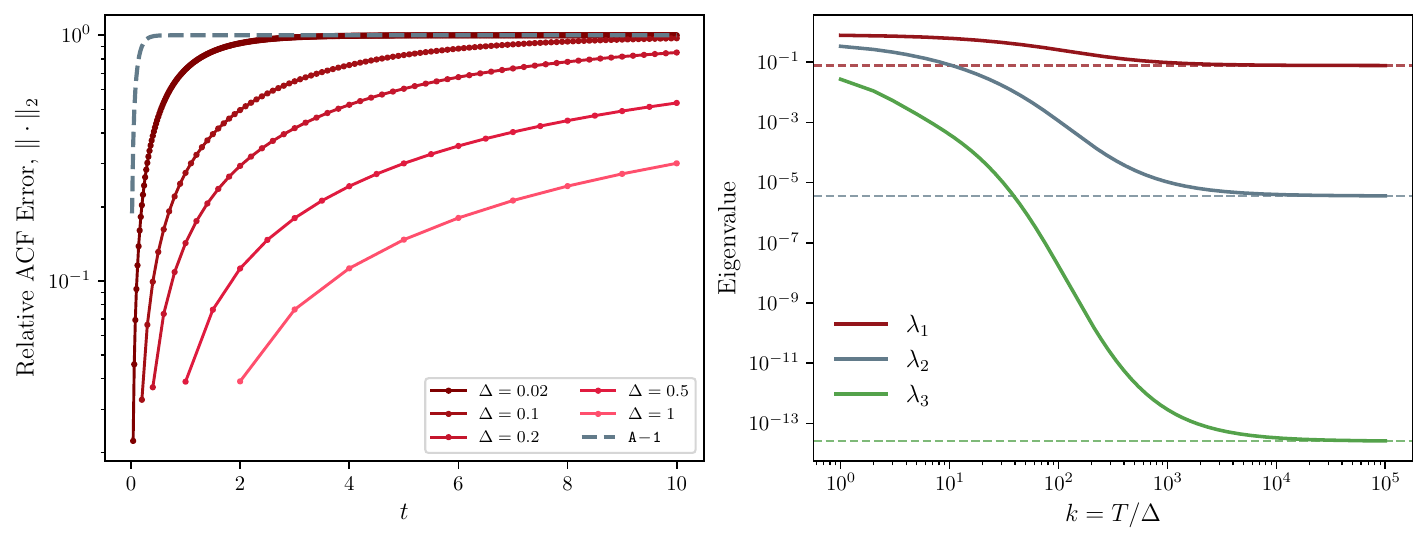}
    \caption{The 40D-to-3D example with $\kappa=1$ and $\ell=500$. Left: relative $2$-norm ACF error for varying $\Delta$. Right: normalized ACF eigenvalues versus step count at $T=0.2$.}
    \label{fig:Ex2_lag}
\end{figure}

In the left panel of Figure~\ref{fig:Ex2_lag}, we fix $\kappa=1$ and $\ell=500$.
We compare $\Delta=0.02,0.1,0.2,0.5,1$ up to $T=10$.
For these parameters, larger lags give smaller errors at the displayed grid
points, while the errors approach the Approach~$\mathtt{A\text{-}1}$ reference
as the lag decreases. This is consistent with the small-lag limit in
Theorem~\ref{thm:A1_limit}.

For the corresponding spectral comparison, we fix $\kappa=1$, $\ell=500$, and $T=0.2$, and vary the number of steps
from $1$ to $10^5$, with $\Delta=T/k$.
The right panel of Figure~\ref{fig:Ex2_lag} shows the three ordered eigenvalues of the
normalized autocovariance in \eqref{Hmc_T}.
At $k=1$, each eigenvalue equals the corresponding eigenvalue of the
normalized full SDE autocovariance. As $k$ increases, the eigenvalues
decrease toward their Approach~$\mathtt{A\text{-}1}$ limits.
These observations agree with Theorem~\ref{thm:monotonicity}.

\begin{figure}[!htbp]
    \centering
    \includegraphics[width=0.99\textwidth]{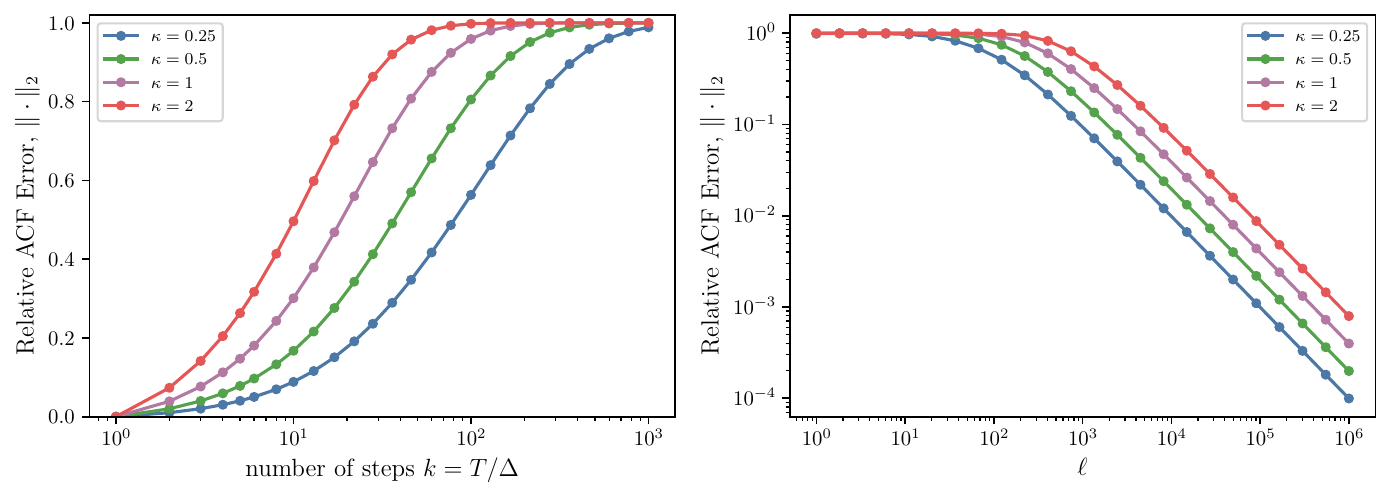}
    \caption{Terminal relative $2$-norm ACF error. Left: varying the step count at $T=5$ and $\ell=500$. Right: varying $\ell$ (separation of time scales) at $T=10$ and $\Delta=0.5$.}
    \label{fig:Ex2_error}
\end{figure}

Figure~\ref{fig:Ex2_error} examines the terminal error for
$\kappa=0.25,0.5,1,2$.
The left panel fixes $T=5$ and $\ell=500$, and varies the number of steps
from $1$ to $10^3$, with $\Delta=T/k$. The error starts at zero, as predicted by
Proposition~\ref{prop:one_step}, and increases as the lag decreases in these
examples. The right panel fixes $T=10$ and $\Delta=0.5$, and varies
$\ell\in[1,10^6]$. The terminal errors decrease for sufficiently large
coupling, showing improved accuracy as the separation of time scales grows.

\subsection{Comparison of Methods}\label{ex:comparison}
\begin{figure}[!htbp]
    \centering
    \includegraphics[width=0.99\textwidth]{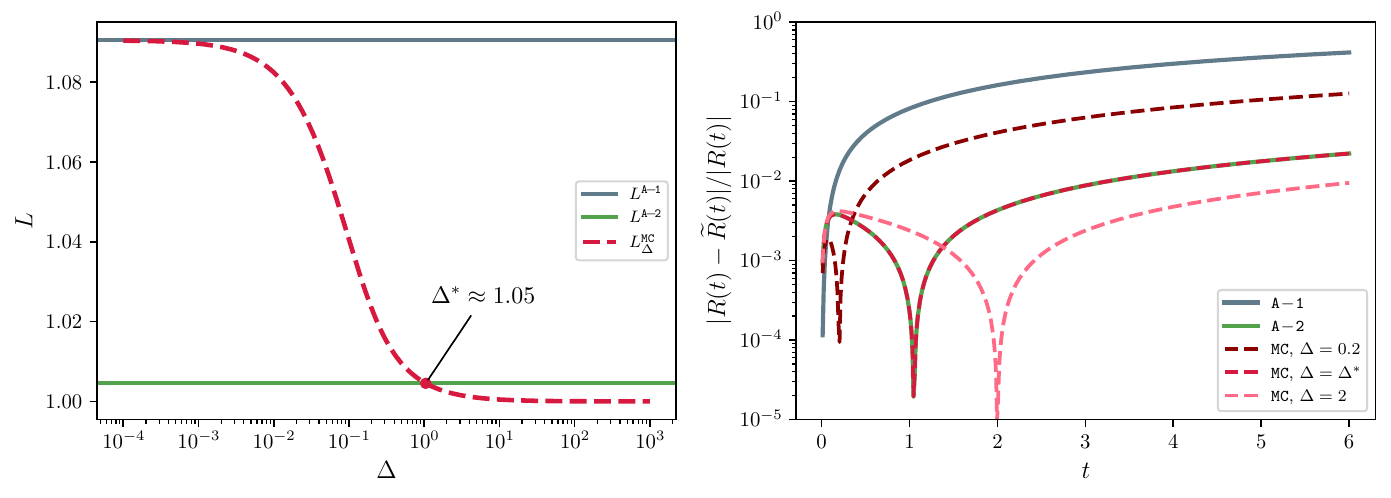}
    \caption{Comparison of the scalar effective models. Left: relaxation rates versus lag. Right: relative autocovariance errors versus time.}
    \label{fig:Ex3}
\end{figure}

We return to the two-dimensional system with $\lambda=20$ and $\theta=0.3$
to compare the effective models from Section~\ref{sec:comparison}.
Figure~\ref{fig:Ex3} uses the relaxation rates in \eqref{Lmc}, \eqref{La1},
and \eqref{La2}, and the autocovariances in \eqref{Hscalar} and \eqref{Hou}.
The left panel varies $\Delta\in[10^{-4},10^3]$.
The embedding rate approaches Approach~$\mathtt{A\text{-}1}$ as the lag
vanishes and matches Approach~$\mathtt{A\text{-}2}$ at
$\Delta^*\approx1.05$, consistent with Lemma~\ref{lem:local} and
Proposition~\ref{prop:scalar_match}.
The right panel compares the relative errors in $R(t)$ up to $T=6$
for $\Delta=0.2,\Delta^*,2$. Since all models have the same equilibrium
variance, these equal the relative errors in $H(t)$.
At the matching lag, the continuous embedding and Approach~$\mathtt{A\text{-}2}$
have coincident error curves. The other lags illustrate how the choice of
lag affects dynamical accuracy.

\FloatBarrier
\section{Conclusion}\label{sec:conclu}
We have studied the stability and statistical properties of a finite-lag
Markov chain for coarse-grained linear SDEs. The chain preserves
the exact equilibrium distribution and matches the autocovariance at the
chosen lag. The analysis quantifies the dynamical error at subsequent times
and clarifies how lag selection affects the approximation. In particular,
the slow--fast analysis shows that lag selection must account for scale
separation, and that a smaller lag need not improve dynamical accuracy.
Through its embedding, the chain is compared with some existing continuous-time
Markovian effective models within a unified framework according to the temporal
statistics they preserve. Numerical examples illustrate these
properties for scalar and vector-valued resolved variables.

\section{Acknowledgments}
The author thanks Prof. Xingjie Helen Li at University of North Carolina at Charlotte for helpful discussions and suggestions.

\section*{AI Statement}
The author used ChatGPT 6 to assist in preparing this manuscript and takes
full responsibility for the accuracy, integrity, and originality of its content.

\appendix
\begingroup
\setlength{\abovedisplayskip}{7pt plus 2pt minus 4pt}
\setlength{\belowdisplayskip}{7pt plus 2pt minus 4pt}
\setlength{\belowdisplayshortskip}{5pt plus 3pt minus 3pt}
\section{Proofs of Statistical Properties}\label{app:statistical}

\begin{proof}[Proof of Lemma~\ref{lem:stability}]
We first prove that the limiting transition law preserves the covariance:
\begin{equation}
    \mathbf{F}_\infty\beta^{-1}\mathbf{B}_0^{-1}\mathbf{F}_\infty^\top+\mathbf{G}_\infty=\beta^{-1}\mathbf{B}_0^{-1}
\end{equation}
Using
\eqref{Ainv} and the block form of $e^{-\mathbf{A}\Delta}$,
\begin{align*}
    [e^{-\mathbf{A}\Delta}\mathbf{A}^{-1}e^{-\mathbf{A}\Delta}]_{\xi\xi}
    &=
    \mathbf{E}_0\mathbf{B}_0^{-1}\mathbf{E}_0^{\top}
    -
    \mathbf{E}_0\mathbf{B}_0^{-1}\boldsymbol{\alpha}
    \mathbf{A}_1^{-1}\mathbf{e}^{\top}
    -
    \mathbf{e}\mathbf{A}_1^{-1}\boldsymbol{\alpha}^{\top}
    \mathbf{B}_0^{-1}\mathbf{E}_0^{\top}
    \\
    &\quad
    +
    \mathbf{e}
    \left(
    \mathbf{A}_1^{-1}
    +
    \mathbf{A}_1^{-1}\boldsymbol{\alpha}^{\top}
    \mathbf{B}_0^{-1}\boldsymbol{\alpha}\mathbf{A}_1^{-1}
    \right)
    \mathbf{e}^{\top} \\
    &=
    \left(\mathbf{E}_0-\mathbf{e}\mathbf{A}_1^{-1}\boldsymbol{\alpha}^{\top}\right)
    \mathbf{B}_0^{-1}
    \left(\mathbf{E}_0-\mathbf{e}\mathbf{A}_1^{-1}\boldsymbol{\alpha}^{\top}\right)^{\top}
    +
    \mathbf{e}\mathbf{A}_1^{-1}\mathbf{e}^{\top} \\
    &=
    \mathbf{F}_\infty\mathbf{B}_0^{-1}\mathbf{F}_\infty^\top
    +
    \mathbf{e}\mathbf{A}_1^{-1}\mathbf{e}^{\top}.
\end{align*}
Since
$\sigma_\Delta=\beta^{-1}
(\mathbf{A}^{-1}-e^{-\mathbf{A}\Delta}\mathbf{A}^{-1}e^{-\mathbf{A}\Delta})$,
the $\xi$-$\xi$ block is
\begin{align*}
    [\sigma_\Delta]_{\xi\xi}
    &=
    \beta^{-1}\mathbf{B}_0^{-1}
    -
    \beta^{-1}\mathbf{F}_\infty\mathbf{B}_0^{-1}\mathbf{F}_\infty^\top
    -
    \beta^{-1}\mathbf{e}\mathbf{A}_1^{-1}\mathbf{e}^{\top}.
\end{align*}
Combining this with the definitions of $\mathbf{F}_\infty$ and
$\mathbf{G}_\infty$ gives
\begin{align*}
    &\mathbf{F}_\infty\beta^{-1}\mathbf{B}_0^{-1}\mathbf{F}_\infty^\top
    +
    \mathbf{G}_\infty \\
    &=
    \mathbf{F}_\infty\beta^{-1}\mathbf{B}_0^{-1}\mathbf{F}_\infty^\top
    +
    [\sigma_\Delta]_{\xi\xi}
    +
    \beta^{-1}\mathbf{e}\mathbf{A}_1^{-1}\mathbf{e}^{\top} \\
    &=
    \mathbf{F}_\infty\beta^{-1}\mathbf{B}_0^{-1}\mathbf{F}_\infty^\top
    +
    \beta^{-1}\mathbf{B}_0^{-1}
    -
    \mathbf{F}_\infty\beta^{-1}\mathbf{B}_0^{-1}\mathbf{F}_\infty^\top
    -
    \beta^{-1}\mathbf{e}\mathbf{A}_1^{-1}\mathbf{e}^{\top}
    +
    \beta^{-1}\mathbf{e}\mathbf{A}_1^{-1}\mathbf{e}^{\top} \\
    &=
    \beta^{-1}\mathbf{B}_0^{-1}.
\end{align*}
For $\Delta>0$, we have
\begin{equation}\label{Gpos}
    \mathbf{G}_\infty
    =
    [\sigma_\Delta]_{\xi\xi}
    +
    \beta^{-1}\mathbf{e}\mathbf{A}_1^{-1}\mathbf{e}^{\top}
    \succ
    \mathbf{0}.
\end{equation}
Here $[\sigma_\Delta]_{\xi\xi}\succ\mathbf{0}$ since it is a principal block of
the positive definite full covariance matrix $\sigma_\Delta$, while
$\beta^{-1}\mathbf{e}\mathbf{A}_1^{-1}\mathbf{e}^{\top}\succeq\mathbf{0}$ since
\(\mathbf{A}_1^{-1}\succ\mathbf{0}\).
Hence
\begin{equation}\label{contraction}
    \mathbf{F}_\infty\beta^{-1}\mathbf{B}_0^{-1}\mathbf{F}_\infty^\top
    =
    \beta^{-1}\mathbf{B}_0^{-1}-\mathbf{G}_\infty
    \prec
    \beta^{-1}\mathbf{B}_0^{-1}.
\end{equation}
To turn this weighted contraction into a statement about eigenvalues, define
\begin{equation}\label{ZF}
    \mathbf{Z}
    :=
    (\beta^{-1}\mathbf{B}_0^{-1})^{-1/2}
    \mathbf{F}_\infty
    (\beta^{-1}\mathbf{B}_0^{-1})^{1/2}.
\end{equation}
Then $\mathbf{Z}$ is similar to $\mathbf{F}_\infty$, and conjugating
\eqref{contraction} by
$(\beta^{-1}\mathbf{B}_0^{-1})^{-1/2}$ gives the expanded inequality
\begin{align}\label{Zbound}
    \mathbf{Z}\mathbf{Z}^\top
    &=
    (\beta^{-1}\mathbf{B}_0^{-1})^{-1/2}
    \mathbf{F}_\infty\beta^{-1}\mathbf{B}_0^{-1}\mathbf{F}_\infty^\top
    (\beta^{-1}\mathbf{B}_0^{-1})^{-1/2} \notag\\
    &\prec
    (\beta^{-1}\mathbf{B}_0^{-1})^{-1/2}
    (\beta^{-1}\mathbf{B}_0^{-1})
    (\beta^{-1}\mathbf{B}_0^{-1})^{-1/2}
    =
    \mathbf{I}_n .
\end{align}
By \eqref{Zbound},
$\sigma_{\max}(\mathbf{Z})^2 = \sigma_{\max}(\mathbf{Z}\mathbf{Z}^{\top}) < \sigma_{\max}(\mathbf{I}_n) = 1$.
Thus by similarity,
\begin{equation}\label{radius}
    \rho(\mathbf{F}_\infty)
    =
    \rho(\mathbf{Z})
    \le
    \|\mathbf{Z}\|_2
    =
    \sigma_{\max}(\mathbf{Z})
    <
    1.
\end{equation}
\end{proof}

\begin{proof}[Proof of Theorem~\ref{thm:int_error}]
We first establish the identity
\begin{equation}\label{int_error_formula}
    \boldsymbol{\mathcal E}_\Delta
    =
    \beta^{-1}\Delta\mathbf{B}_0^{-1/2}
    \left[
    \mathbf{U}^\top
    (\mathbf{I}_N-e^{-\mathbf{A}\Delta})^{-1}
    \mathbf{U}
    -
    (\mathbf{I}_n-\mathbf{H}(\Delta))^{-1}
    \right]
    \mathbf{B}_0^{-1/2}.
\end{equation}
\RV{Equations \eqref{U} and
\eqref{H} give
\begin{align}
    \beta\mathbf{B}_0^{1/2}
    \boldsymbol{\mathcal E}_\Delta
    \mathbf{B}_0^{1/2}
    &=
    \Delta\sum_{k=0}^{\infty}
    \left[
    \mathbf{U}^\top e^{-\mathbf{A}k\Delta}\mathbf{U}
    -
    \mathbf{H}(\Delta)^k
    \right] \notag\\
    &=
    \Delta
    \left[
    \mathbf{U}^\top
    (\mathbf{I}_N-e^{-\mathbf{A}\Delta})^{-1}
    \mathbf{U}
    -
    (\mathbf{I}_n-\mathbf{H}(\Delta))^{-1}
    \right].\label{int_sum}
\end{align}
Since $\mathbf{A}\succ\mathbf{0}$ and $\Delta>0$,
$\mathbf{I}_N-e^{-\mathbf{A}\Delta}\succ\mathbf{0}$. Its block representation in the orthogonal basis $(\mathbf{U},\mathbf{V})$ is
\begin{equation}\label{int_block}
    \begin{pmatrix}\mathbf{U}^\top\\ \mathbf{V}^\top\end{pmatrix}
    (\mathbf{I}_N-e^{-\mathbf{A}\Delta})
    \begin{pmatrix}\mathbf{U}&\mathbf{V}\end{pmatrix}
    =
    \begin{pmatrix}
        \mathbf{I}_n-\mathbf{H}(\Delta)
        &-\mathbf{U}^\top e^{-\mathbf{A}\Delta}\mathbf{V}\\
        -\mathbf{V}^\top e^{-\mathbf{A}\Delta}\mathbf{U}
        &\mathbf{I}_{N-n}-\mathbf{V}^\top e^{-\mathbf{A}\Delta}\mathbf{V}
    \end{pmatrix}
    \succ\mathbf{0}.
\end{equation}
Its lower-right block is positive definite. Then we consider the upper-left block of inverse of \eqref{int_block} using block-inverse formula
\begin{equation}\label{int_Schur}
    \begin{aligned}
        \mathbf{U}^\top
        (\mathbf{I}_N-e^{-\mathbf{A}\Delta})^{-1}
        \mathbf{U}
        &=
        \Bigl[
        \mathbf{I}_n-\mathbf{H}(\Delta)
        -\mathbf{U}^\top e^{-\mathbf{A}\Delta}\mathbf{V}\\
        &\qquad
        (\mathbf{I}_{N-n}-\mathbf{V}^\top e^{-\mathbf{A}\Delta}\mathbf{V})^{-1}
        \mathbf{V}^\top e^{-\mathbf{A}\Delta}\mathbf{U}
        \Bigr]^{-1}\succeq
        (\mathbf{I}_n-\mathbf{H}(\Delta))^{-1}.
    \end{aligned}
\end{equation}
}
Indeed, the cross term subtracted inside the inverse is positive semidefinite,
and inversion reverses the Loewner order on positive definite matrices.
Combining \eqref{int_sum} and
\eqref{int_Schur} proves \eqref{int_error_nonnegative} and
shows that repeated Markov transitions systematically underestimate the
integrated autocovariance.
\end{proof}

\begin{proof}[Proof of Corollary~\ref{cor:int_bound}]
\RV{Apply the resolvent identity
\begin{equation}\label{resolvent}
    \mathbf{X}^{-1}-\mathbf{Y}^{-1}
    =\mathbf{X}^{-1}(\mathbf{Y}-\mathbf{X})\mathbf{Y}^{-1}
\end{equation}
to the two inverses in \eqref{int_Schur}. The full resolvent,
the cross correction, and the compressed resolvent satisfy, respectively,
\begin{align}
    \|\mathbf{U}^\top(\mathbf{I}_N-e^{-\mathbf{A}\Delta})^{-1}\mathbf{U}\|_2
    &\le\frac{1}{1-r_\Delta},\notag\\
    \|\mathbf{U}^\top e^{-\mathbf{A}\Delta}\mathbf{V}
    (\mathbf{I}_{N-n}-\mathbf{V}^\top e^{-\mathbf{A}\Delta}\mathbf{V})^{-1}
    \mathbf{V}^\top e^{-\mathbf{A}\Delta}\mathbf{U}\|_2
    &\le\frac{\eta_\Delta^2}{1-\ell_\Delta},\notag\\
    \|(\mathbf{I}_n-\mathbf{H}(\Delta))^{-1}\|_2
    &=\frac{1}{1-h_\Delta}.
    \label{resolvent_bounds}
\end{align}
Therefore, combining \eqref{int_sum}, \eqref{int_Schur}, \eqref{resolvent} and \eqref{resolvent_bounds} gives
\begin{equation}\label{int_norm}
    \|\mathbf{U}^\top(\mathbf{I}_N-e^{-\mathbf{A}\Delta})^{-1}\mathbf{U}
    -(\mathbf{I}_n-\mathbf{H}(\Delta))^{-1}\|_2
    \le\frac{\eta_\Delta^2}{(1-r_\Delta)(1-\ell_\Delta)(1-h_\Delta)}.
\end{equation}
Substitution into \eqref{int_error_formula} proves
\eqref{int_bound}. Moreover,
\begin{equation}\label{decay_bounds}
    r_\Delta=e^{-\lambda_{\min}(\mathbf{A})\Delta},
    \qquad
    h_\Delta\le r_\Delta,
    \qquad
    \ell_\Delta\le r_\Delta,
    \qquad
    \eta_\Delta^2 \le h_\Delta \ell_\Delta \le r_\Delta^2.
\end{equation}
Substitution into \eqref{int_bound} gives
\eqref{int_decay}, whose right-hand side converges
to zero as $\Delta\to\infty$.}
\end{proof}

\begin{proof}[Proof of Theorem~\ref{thm:monotonicity}]
Using the definition \eqref{U}, we have
\begin{equation}\label{Hpower}
    \mathbf{H}(\Delta_k)
    =
    \mathbf{U}^\top e^{-\mathbf{A}T/k}\mathbf{U}.
\end{equation}
Let $f_k(x)=x^{k/(k+1)}$ and
$\Phi(\mathbf{M})=\mathbf{U}^\top\mathbf{M}\mathbf{U}$. The function $f_k$ is
operator concave on $(0,\infty)$, and $\Phi$ is a unital positive map since
$\mathbf{U}^\top\mathbf{U}=\mathbf{I}_n$. Applying operator Jensen's
inequality~\cite[Theorem~2.1]{hansen2003jensen},
$\Phi(f_k(\mathbf{M}))\preceq f_k(\Phi(\mathbf{M}))$, to
$\mathbf{M}=e^{-\mathbf{A}T/k}$ and using
\eqref{Hpower} give
\begin{align}\label{Hconcavity}
    \mathbf{H}(\Delta_{k+1})
    =
    \mathbf{U}^\top
    (e^{-\mathbf{A}T/k})^{k/(k+1)}
    \mathbf{U}
    \preceq
    \left(
    \mathbf{U}^\top e^{-\mathbf{A}T/k}\mathbf{U}
    \right)^{k/(k+1)}
    =
    \mathbf{H}(\Delta_k)^{k/(k+1)}.
\end{align}
\RV{By eigenvalue monotonicity under the Loewner order,
\eqref{Hconcavity} gives
\[
    \lambda_j\!\left(\mathbf H(\Delta_{k+1})\right)
    \le
    \lambda_j\!\left(\mathbf H(\Delta_k)\right)^{k/(k+1)}.
\]
Raising this scalar inequality to the power $k+1$ and using
\eqref{Hmc_T} yield}
\begin{equation}\label{eig_order}
    \lambda_j(\mathbf{H}(\Delta_{k+1})^{k+1})
    \le
    \lambda_j(\mathbf{H}(\Delta_k)^k),
\end{equation}
which proves \eqref{monotonicity}. Finally,
\begin{equation}\label{Hexpand}
    \frac{\mathbf{I}_n-\mathbf{H}(t)}{t}=\mathbf{U}^\top\frac{\RV{\mathbf{I}_N}-e^{-\mathbf{A}t}}{t}\mathbf{U}
    \longrightarrow
    \mathbf{U}^\top\mathbf{A}\mathbf{U}
    =
    \mathbf{B}_0,
    \qquad t\to0.
\end{equation}
Taking $t=T/k$ in \eqref{Hexpand}, we obtain $k\left(\mathbf{I}_n-\mathbf{H}(T/k)\right)\longrightarrow T\mathbf{B}_0$. \RV{The standard matrix exponential limit then gives}
\begin{align}
    \mathbf{H}(\Delta_k)^k
    =
    \left[
    \mathbf{I}_n-
    \frac{k\left(\mathbf{I}_n-\mathbf{H}(T/k)\right)}{k}
    \right]^k
    \longrightarrow e^{-\mathbf{B}_0T},
    \label{Hlimit}
\end{align}
which proves \eqref{Hmc_limit}.
\end{proof}

\section{Proofs for the Two-Dimensional Slow--Fast System}\label{app:slow_fast}

\RV{We collect the exact identities for the system in \eqref{model2d}. For
a fixed lag time $\Delta>0$, define}
\begin{equation}\label{param2d}
    \RV{
    \rho=\frac{s^2}{c^2}=\tan^2\theta,\qquad
    b=\frac{\rho}{\lambda},\qquad
    r=\frac{q(\Delta)}{p(\Delta)}=e^{-(\lambda-1)\Delta}.}
\end{equation}
\RV{Then the one-step coefficient can be written equivalently as}
\begin{equation}\label{F2d}
    \RV{F_\infty(\Delta,\lambda)
    = \frac{c^2\lambda p(\Delta) + s^2 q(\Delta)}{s^2 + c^2\lambda}
    = \frac{p(\Delta)+bq(\Delta)}{1+b}
    =p(\Delta)\frac{1+br}{1+b}.}
\end{equation}
\RV{The effective equilibrium autocovariance is}
\begin{equation}\label{Rmc2d}
    \RV{\widetilde{R}(k\Delta)
    =\beta^{-1}c^2p(\Delta)^k(1+b)
    \left(\frac{1+br}{1+b}\right)^k,}
\end{equation}
\RV{whereas the exact autocovariance is}
\begin{equation} \label{exactR}
    \RV{R(k\Delta)
    =\beta^{-1}c^2p(\Delta)^k(1+br^k).}
\end{equation}
\RV{In particular, their exact ratio is}
\begin{equation}\label{Rratio}
    \RV{
    \frac{\widetilde{R}(k\Delta)}{R(k\Delta)}
    =F_\infty(\Delta,\lambda)^k e^{k\Delta}
    \frac{1+b}{1+br^k}
    =\frac{1+b}{1+br^k}
    \left(\frac{1+br}{1+b}\right)^k.}
\end{equation}

\begin{proof}[Proof of Theorem~\ref{thm:fixed_lag}]
\RV{For fixed $\Delta>0$ and $k$, the quantity
$r=e^{-(\lambda-1)\Delta}$ is exponentially small as $\lambda\to\infty$,
whereas $b=\rho\lambda^{-1}$. Expanding the exact ratio
\eqref{Rratio} factor by factor gives}
\begin{equation}\label{ratio2d_fixed}
    \RV{
    \begin{aligned}
    \frac{1+b}{1+br^k}
    &=(1+b)\bigl(1-br^k+O(b^2r^{2k})\bigr)
      =1+b+O(\lambda^{-2}),\\
    \left(\frac{1+br}{1+b}\right)^k
    &=\bigl(1+b(r-1)+O(b^2)\bigr)^k
      =1-kb+O(\lambda^{-2}),\\
    \frac{\widetilde{R}(k\Delta)}{R(k\Delta)}
    &=\bigl(1+b+O(\lambda^{-2})\bigr)
      \bigl(1-kb+O(\lambda^{-2})\bigr)\\
    &=1+\frac{(1-k)\rho}{\lambda}+O(\lambda^{-2}).
    \end{aligned}}
\end{equation}
\RV{Since $R(k\Delta)>0$ and $\rho=\tan^2\theta$, taking the absolute
difference from $1$ proves the result.}
\end{proof}

\begin{proof}[Proof of Theorem~\ref{thm:small_lag}]
\RV{We first record the expansion common to all three cases. From
\eqref{F2d},}
\begin{equation}\label{rate2d}
    \RV{
    \frac{1-F_\infty(\Delta,\lambda)}{\Delta}
    =
    \frac{
    \dfrac{1-e^{-\Delta}}{\Delta}
    +\rho\dfrac{1-e^{-\lambda\Delta}}{\lambda\Delta}}
    {1+\rho\lambda^{-1}}.}
\end{equation}
\RV{The function $(1-e^{-x})/x$ is uniformly bounded, so
$0 \leq 1-F_\infty(\Delta,\lambda)\leq (1+\rho)\Delta$. Since $k\Delta\to\tau$, the
Taylor expansion of the logarithm gives}
\begin{equation}\label{logF2d}
    \RV{
    \begin{aligned}
    k\log F_\infty
    &=-k(1-F_\infty)+O\!\left(k(1-F_\infty)^2\right)=-k\Delta\frac{1-F_\infty}{\Delta}+O(\Delta).
    \end{aligned}}
\end{equation}
\RV{The factor multiplying $F_\infty^k e^{k\Delta}$ in
\eqref{Rratio} lies between $1$ and $1+b$, so it is
$1+O(\lambda^{-1})$, uniformly in $k$ and $\Delta$.}

\RV{\textit{Case 1: $\lambda\Delta\to0$.}
Expanding the two exponentials in
\eqref{rate2d} gives}
\begin{equation}\label{rate2d_small}
    \RV{
    \frac{1-F_\infty}{\Delta}
    =1+\rho+O(\lambda\Delta+\lambda^{-1}).}
\end{equation}
\RV{Substituting \eqref{rate2d_small} into \eqref{logF2d}, and then using
\eqref{Rratio}, yields}
\begin{equation}\label{ratio2d_small}
    \RV{
    \frac{\widetilde{R}(k\Delta)}{R(k\Delta)}
    =e^{-\rho\tau}
    \left[1+O\!\left(|k\Delta-\tau|
    +\lambda\Delta+\lambda^{-1}\right)\right].}
\end{equation}
\RV{Taking the absolute difference from $1$ proves the first case.}

\RV{\textit{Case 2: $\lambda\Delta\to\ell\in(0,\infty)$.}
Since the function
$x\mapsto(1-e^{-x})/x$ is smooth near $\ell$, and
$\lambda^{-1}=O(\Delta)$ in this regime, \eqref{rate2d} gives}
\begin{equation}\label{rate2d_mid}
    \RV{
    \frac{1-F_\infty}{\Delta}
    =1+\rho\frac{1-e^{-\ell}}{\ell}
    +O\!\left(|\lambda\Delta-\ell|+\Delta\right).}
\end{equation}
\RV{Equations \eqref{logF2d} and \eqref{Rratio} therefore imply}
\begin{equation}\label{ratio2d_mid}
    \RV{
    \frac{\widetilde{R}(k\Delta)}{R(k\Delta)}
    =e^{-\rho\tau(1-e^{-\ell})/\ell}
    \left[1+O\!\left(|k\Delta-\tau|
    +|\lambda\Delta-\ell|+\Delta\right)\right],}
\end{equation}
\RV{which proves the second case after taking the absolute difference from
$1$.}

\RV{\textit{Case 3: $\lambda\Delta\to\infty$.}
To identify the leading term, we use the exact factorization rather than
the Taylor expansion \eqref{logF2d}. Taking logarithms in
\eqref{Rratio} yields}
\begin{equation}\label{logratio2d}
    \RV{
    \begin{aligned}
    \log\frac{\widetilde{R}(k\Delta)}{R(k\Delta)}
    &=k\log(1+br)-(k-1)\log(1+b)-\log(1+br^k)\\
    &=-\frac{\rho\tau}{\lambda\Delta}
    +O\!\left(
    \frac{|k\Delta-\tau|}{\lambda\Delta}
    +\lambda^{-1}+\frac{1}{(\lambda\Delta)^2}
    \right).
    \end{aligned}}
\end{equation}
\RV{Here we used $kbr=O(e^{-\lambda\Delta}/(\lambda\Delta))$,
$kb^2=O(\lambda^{-1})$, and $br^k=O(\lambda^{-1})$.
Exponentiating \eqref{logratio2d} and using $1-e^{-x}=x+O(x^2)$ gives}
\begin{equation}\label{error2d_large}
    \RV{
    \frac{|\widetilde{R}(k\Delta)-R(k\Delta)|}{R(k\Delta)}
    =\frac{\rho\tau}{\lambda\Delta}
    +O\!\left(
    \frac{|k\Delta-\tau|}{\lambda\Delta}
    +\lambda^{-1}+\frac{1}{(\lambda\Delta)^2}
    \right).}
\end{equation}
\RV{Since each term in the remainder is
$o((\lambda\Delta)^{-1})$ and $\rho=\tan^2\theta$, this proves the third
case with the stated leading term.}
\end{proof}

\begin{proof}[Proof of Theorem~\ref{thm:finite_time}]
\RV{The sign of the error follows from the convexity of $x\mapsto x^k$
on $[0,1]$. Indeed, by
\eqref{F2d}--\eqref{exactR},}
\begin{equation}\label{Rorder}
    \RV{
    \widetilde{R}(k\Delta)
    =\beta^{-1}c^2(1+b)
    \left(\frac{p(\Delta)+bq(\Delta)}{1+b}\right)^k
    \le \beta^{-1}c^2\left(p(\Delta)^k+bq(\Delta)^k\right)
    =R(k\Delta).}
\end{equation}
\RV{Hence $0\le R(k\Delta)-\widetilde{R}(k\Delta)$ for every $k\ge1$.}

\RV{For the upper bound, $0<r=e^{-(\lambda-1)\Delta}<1$ because $\lambda>1$. Thus the common
ratio \eqref{Rratio} gives, for $t=k\Delta$,}
\begin{equation}\label{Rratio_bound}
    \RV{
    \frac{\widetilde{R}(t)}{R(t)}
    \ge \left(\frac{1+br}{1+b}\right)^k.}
\end{equation}
\RV{Moreover,}
\begin{equation}\label{log_bound}
    \RV{
    0\le \log(1+b)-\log(1+br)
    =\int_{br}^{b}\frac{1}{1+x}\,\diff x
    \le b(1-r).}
\end{equation}
\RV{Therefore, since $t=k\Delta\le T$,}
\begin{equation}\label{Rratio_exp}
    \RV{
    \frac{\widetilde{R}(t)}{R(t)}
    \ge e^{-kb(1-r)}
    \ge
    e^{-T\rho(1-e^{-(\lambda-1)\Delta})/(\lambda\Delta)}.}
\end{equation}
\RV{Combining this estimate with \eqref{Rorder}, taking the
supremum over $1\le k\le\lfloor T/\Delta\rfloor$, and using
$\rho=\tan^2\theta$ proves \eqref{error2d_bound}.}
\end{proof}

\section{Proof of Lemma~\ref{lem:local}}\label{app:local}

\begin{proof}
\RV{By \eqref{Hexpand},
\begin{equation}\label{Hsmall}
    \frac{\mathbf{H}(\Delta)-\mathbf{I}_n}{\Delta}
    \longrightarrow
    -\mathbf{B}_0.
\end{equation}
In particular, $\|\mathbf{H}(\Delta)-\mathbf{I}_n\|_2\to0$. Hence, for
sufficiently small $\Delta$, the matrix Taylor series for the principal
logarithm converges and gives
\begin{align}\label{logH}
    \log\mathbf{H}(\Delta)
    &=
    \sum_{m=1}^{\infty}
    \frac{(-1)^{m+1}}{m}
    \bigl(\mathbf{H}(\Delta)-\mathbf{I}_n\bigr)^m \notag\\
    &=
    \bigl(\mathbf{H}(\Delta)-\mathbf{I}_n\bigr)
    \left[
    \mathbf{I}_n+
    \sum_{m=2}^{\infty}
    \frac{(-1)^{m+1}}{m}
    \bigl(\mathbf{H}(\Delta)-\mathbf{I}_n\bigr)^{m-1}
    \right].
\end{align}
Moreover,
\begin{align}\label{logH_bound}
    &\left\|
    \sum_{m=2}^{\infty}
    \frac{(-1)^{m+1}}{m}
    \bigl(\mathbf{H}(\Delta)-\mathbf{I}_n\bigr)^{m-1}
    \right\|_2 \notag\\
    &\qquad\le
    \sum_{m=2}^{\infty}
    \|\mathbf{H}(\Delta)-\mathbf{I}_n\|_2^{m-1}
    =
    \frac{\|\mathbf{H}(\Delta)-\mathbf{I}_n\|_2}
    {1-\|\mathbf{H}(\Delta)-\mathbf{I}_n\|_2}
    \longrightarrow0.
\end{align}
Dividing \eqref{logH} by $\Delta$ and using
\eqref{Hsmall} therefore yields
\begin{equation}\label{logH_limit}
    \frac{\log\mathbf{H}(\Delta)}{\Delta}
    \longrightarrow-\mathbf{B}_0.
\end{equation}
Thus we proved
\eqref{Lmc_limit}.}
\end{proof}

\endgroup

\bibliographystyle{siamplain}
\bibliography{a}

@article{hudson2024dynamical,
  title={Dynamical properties of coarse-grained linear SDEs},
  author={Hudson, Thomas and Li, Xingjie Helen},
  journal={Multiscale Modeling \& Simulation},
  volume={22},
  number={1},
  pages={406--435},
  year={2024},
  publisher={SIAM}
}

@article{lei2016data,
  title={Data-driven parameterization of the generalized Langevin equation},
  author={Lei, Huan and Baker, Nathan A and Li, Xiantao},
  journal={Proceedings of the National Academy of Sciences},
  volume={113},
  number={50},
  pages={14183--14188},
  year={2016},
  publisher={National Academy of Sciences}
}

@article{chen2024datadriven,
  author    = {Yuan Chen and Dongbin Xiu},
  title     = {Data-Driven Effective Modeling of Multiscale Stochastic Dynamical Systems},
  journal   = {Applied Mathematics for Modern Challenges},
  volume    = {2},
  number    = {4},
  pages     = {348--364},
  year      = {2024},
  doi       = {10.3934/ammc.2024018},
}

@article{zwanzig1973nonlinear,
  title={Nonlinear generalized Langevin equations},
  author={Zwanzig, Robert},
  journal={Journal of Statistical Physics},
  volume={9},
  number={3},
  pages={215--220},
  year={1973},
  publisher={Springer}
}

@book{zwanzig2001nonequilibrium,
  title={Nonequilibrium statistical mechanics},
  author={Zwanzig, Robert},
  year={2001},
  publisher={Oxford university press}
}

@article{legoll2019effective,
  title={Effective dynamics for non-reversible stochastic differential equations: a quantitative study},
  author={Legoll, Fr{\'e}d{\'e}ric and Leli{\`e}vre, Tony and Sharma, Upanshu},
  journal={Nonlinearity},
  volume={32},
  number={12},
  pages={4779--4816},
  year={2019},
  publisher={IOP Publishing}
}

@article{legoll2010effective,
  title={Effective dynamics using conditional expectations},
  author={Legoll, Fr{\'e}d{\'e}ric and Lelievre, Tony},
  journal={Nonlinearity},
  volume={23},
  number={9},
  pages={2131--2163},
  year={2010}
}

@article{nuske2021spectral,
  title={Spectral Properties of Effective Dynamics from Conditional Expectations},
  author={N{\"u}ske, Feliks and Koltai, P{\'e}ter and Boninsegna, Lorenzo and Clementi, Cecilia},
  journal={Entropy},
  volume={23},
  number={2},
  pages={134},
  year={2021},
  publisher={MDPI}
}

@book{pavliotis2008multiscale,
  title={Multiscale Methods: Averaging and Homogenization},
  author={Pavliotis, Grigorios A. and Stuart, Andrew M.},
  series={Texts in Applied Mathematics},
  volume={53},
  year={2008},
  publisher={Springer}
}

@book{kloeden1992numerical,
  title={Numerical Solution of Stochastic Differential Equations},
  author={Kloeden, Peter E. and Platen, Eckhard},
  series={Applications of Mathematics},
  volume={23},
  year={1992},
  publisher={Springer-Verlag},
  address={Berlin},
  doi={10.1007/978-3-662-12616-5}
}

@article{noid2013perspective,
  author={Noid, W. G.},
  title={Perspective: Coarse-grained models for biomolecular systems},
  journal={The Journal of Chemical Physics},
  volume={139}, number={9}, pages={090901}, year={2013},
  doi={10.1063/1.4818908}
}

@article{best2010coordinate,
  author={Best, Robert B. and Hummer, Gerhard},
  title={Coordinate-dependent diffusion in protein folding},
  journal={Proceedings of the National Academy of Sciences},
  volume={107}, number={3}, pages={1088--1093}, year={2010},
  doi={10.1073/pnas.0910390107}
}

@article{prinz2011markov,
  author={Prinz, Jan-Hendrik and Wu, Hao and Sarich, Marco and Keller, Bettina and Senne, Martin and Held, Martin and Chodera, John D. and Sch{\"u}tte, Christof and No{\'e}, Frank},
  title={{Markov} models of molecular kinetics: Generation and validation},
  journal={The Journal of Chemical Physics},
  volume={134}, number={17}, pages={174105}, year={2011},
  doi={10.1063/1.3565032}
}

@article{majda2003systematic,
  author={Majda, Andrew J. and Timofeyev, Ilya and Vanden-Eijnden, Eric},
  title={Systematic strategies for stochastic mode reduction in climate},
  journal={Journal of the Atmospheric Sciences},
  volume={60}, number={14}, pages={1705--1722}, year={2003},
  doi={10.1175/1520-0469(2003)060<1705:SSFSMR>2.0.CO;2}
}

@article{mori1965transport,
  author={Mori, Hazime},
  title={Transport, collective motion, and {Brownian} motion},
  journal={Progress of Theoretical Physics},
  volume={33}, number={3}, pages={423--455}, year={1965},
  doi={10.1143/ptp.33.423}
}

@article{chorin2000optimal,
  author={Chorin, Alexandre J. and Hald, Ole H. and Kupferman, Raz},
  title={Optimal prediction and the {Mori--Zwanzig} representation of irreversible processes},
  journal={Proceedings of the National Academy of Sciences},
  volume={97}, number={7}, pages={2968--2973}, year={2000},
  doi={10.1073/pnas.97.7.2968}
}

@article{jung2018generalized,
  author={Jung, Gerhard and Hanke, Martin and Schmid, Friederike},
  title={Generalized {Langevin} dynamics: Construction and numerical integration of non-{Markovian} particle-based models},
  journal={Soft Matter},
  volume={14}, number={46}, pages={9368--9382}, year={2018},
  doi={10.1039/C8SM01817K}
}

@article{grogan2020data,
  author={Grogan, Francesca and Lei, Huan and Li, Xiantao and Baker, Nathan A.},
  title={Data-driven molecular modeling with the generalized {Langevin} equation},
  journal={Journal of Computational Physics},
  volume={418}, pages={109633}, year={2020},
  doi={10.1016/j.jcp.2020.109633}
}

@article{ge2024data,
  author={Ge, Pei and Zhang, Zhongqiang and Lei, Huan},
  title={Data-driven learning of the generalized {Langevin} equation with state-dependent memory},
  journal={Physical Review Letters},
  volume={133}, number={7}, pages={077301}, year={2024},
  doi={10.1103/PhysRevLett.133.077301}
}

@article{klippenstein2021introducing,
  author={Klippenstein, Viktor and Tripathy, Madhusmita and Jung, Gerhard and Schmid, Friederike and van der Vegt, Nico F. A.},
  title={Introducing memory in coarse-grained molecular simulations},
  journal={The Journal of Physical Chemistry B},
  volume={125}, number={19}, pages={4931--4954}, year={2021},
  doi={10.1021/acs.jpcb.1c01120}
}

@article{majda2001mathematical,
  author={Majda, Andrew J. and Timofeyev, Ilya and Vanden-Eijnden, Eric},
  title={A mathematical framework for stochastic climate models},
  journal={Communications on Pure and Applied Mathematics},
  volume={54}, number={8}, pages={891--974}, year={2001},
  doi={10.1002/cpa.1014}
}

@article{boninsegna2018sparse,
  author={Boninsegna, Lorenzo and N{\"u}ske, Feliks and Clementi, Cecilia},
  title={Sparse learning of stochastic dynamical equations},
  journal={The Journal of Chemical Physics},
  volume={148}, number={24}, pages={241723}, year={2018},
  doi={10.1063/1.5018409}
}

@article{schutte1999direct,
  author={Sch{\"u}tte, Christof and Fischer, Alexander and Huisinga, Wilhelm and Deuflhard, Peter},
  title={A direct approach to conformational dynamics based on hybrid {Monte Carlo}},
  journal={Journal of Computational Physics},
  volume={151}, number={1}, pages={146--168}, year={1999},
  doi={10.1006/jcph.1999.6231}
}

@article{darve2009computing,
  author={Darve, Eric and Solomon, Jose and Kia, Amirali},
  title={Computing generalized {Langevin} equations and generalized {Fokker--Planck} equations},
  journal={Proceedings of the National Academy of Sciences},
  volume={106}, number={27}, pages={10884--10889}, year={2009},
  doi={10.1073/pnas.0902633106}
}

@article{hijon2010mori,
  author={Hij{\'o}n, Carmen and Espa{\~n}ol, Pep and Vanden-Eijnden, Eric and Delgado-Buscalioni, Rafael},
  title={{Mori--Zwanzig} formalism as a practical computational tool},
  journal={Faraday Discussions},
  volume={144}, pages={301--322}, year={2010},
  doi={10.1039/B902479B}
}

@article{hudson2020coarse,
  author={Hudson, Thomas and Li, Xingjie Helen},
  title={Coarse-graining of overdamped {Langevin} dynamics via the {Mori--Zwanzig} formalism},
  journal={Multiscale Modeling \& Simulation},
  volume={18}, number={2}, pages={1113--1135}, year={2020},
  doi={10.1137/18M1222533}
}

@article{vroylandt2022position,
  author={Vroylandt, Hadrien and Monmarch{\'e}, Pierre},
  title={Position-dependent memory kernel in generalized {Langevin} equations: Theory and numerical estimation},
  journal={The Journal of Chemical Physics},
  volume={156}, number={24}, pages={244105}, year={2022},
  doi={10.1063/5.0094566}
}

@article{kevrekidis2003equation,
  author={Kevrekidis, Ioannis G. and Gear, C. William and Hyman, James M. and Kevrekidis, Panagiotis G. and Runborg, Olof and Theodoropoulos, Constantinos},
  title={Equation-free, coarse-grained multiscale computation: Enabling microscopic simulators to perform system-level analysis},
  journal={Communications in Mathematical Sciences},
  volume={1}, number={4}, pages={715--762}, year={2003},
  doi={10.4310/CMS.2003.v1.n4.a5}
}

@article{e2007heterogeneous,
  author={E, Weinan and Engquist, Bjorn and Li, Xiantao and Ren, Weiqing and Vanden-Eijnden, Eric},
  title={Heterogeneous multiscale methods: A review},
  journal={Communications in Computational Physics},
  volume={2}, number={3}, pages={367--450}, year={2007},
  doi={10.4208/cicp.2007.v2.p367}
}

@article{legoll2017pathwise,
  author={Legoll, Fr{\'e}d{\'e}ric and Leli{\`e}vre, Tony and Olla, Stefano},
  title={Pathwise estimates for an effective dynamics},
  journal={Stochastic Processes and their Applications},
  volume={127}, number={9}, pages={2841--2863}, year={2017},
  doi={10.1016/j.spa.2017.01.001}
}

@article{lelievre2019pathwise,
  author={Leli{\`e}vre, Tony and Zhang, Wei},
  title={Pathwise estimates for effective dynamics: The case of nonlinear vectorial reaction coordinates},
  journal={Multiscale Modeling \& Simulation},
  volume={17}, number={3}, pages={1019--1051}, year={2019},
  doi={10.1137/18M1186034}
}

@article{harmandaris2016path,
  author={Harmandaris, Vagelis and Kalligiannaki, Evangelia and Katsoulakis, Markos and Plech{\'a}{\v c}, Petr},
  title={Path-space variational inference for non-equilibrium coarse-grained systems},
  journal={Journal of Computational Physics},
  volume={314}, pages={355--383}, year={2016},
  doi={10.1016/j.jcp.2016.03.021}
}

@article{sachs2025machine,
  author={Sachs, Matthias and Stark, Wojciech G. and Maurer, Reinhard J. and Ortner, Christoph},
  title={Machine learning configuration-dependent friction tensors in {Langevin} heatbaths},
  journal={Machine Learning: Science and Technology},
  volume={6}, number={1}, pages={015016}, year={2025}
}

@article{zhang2016effective,
  author={Zhang, Wei and Hartmann, Carsten and Sch{\"u}tte, Christof},
  title={Effective dynamics along given reaction coordinates, and reaction rate theory},
  journal={Faraday Discussions},
  volume={195}, pages={365--394}, year={2016},
  doi={10.1039/C6FD00147E}
}

@article{klus2018data,
  author={Klus, Stefan and N{\"u}ske, Feliks and Koltai, P{\'e}ter and Wu, Hao and Kevrekidis, Ioannis and Sch{\"u}tte, Christof and No{\'e}, Frank},
  title={Data-driven model reduction and transfer operator approximation},
  journal={Journal of Nonlinear Science},
  volume={28}, pages={985--1010}, year={2018},
  doi={10.1007/s00332-017-9437-7}
}

@article{zhang2025finding,
  author={Zhang, Wei and Sch{\"u}tte, Christof},
  title={On finding optimal collective variables for complex systems by minimizing the deviation between effective and full dynamics},
  journal={Multiscale Modeling \& Simulation},
  volume={23}, number={2}, pages={924--958}, year={2025},
  doi={10.1137/24M1658917}
}

@article{perez2013identification,
  author={P{\'e}rez-Hern{\'a}ndez, Guillermo and Paul, Fabian and Giorgino, Toni and de Fabritiis, Gianni and No{\'e}, Frank},
  title={Identification of slow molecular order parameters for {Markov} model construction},
  journal={The Journal of Chemical Physics},
  volume={139}, number={1}, pages={015102}, year={2013}
}

@article{noe2013variational,
  author={No{\'e}, Frank and N{\"u}ske, Feliks},
  title={A variational approach to modeling slow processes in stochastic dynamical systems},
  journal={Multiscale Modeling \& Simulation},
  volume={11}, number={2}, pages={635--655}, year={2013}
}

@article{chekroun2020ruelle,
  author={Chekroun, Micka{\"e}l D. and Tantet, Alexis and Dijkstra, Henk A. and Neelin, J. David},
  title={{Ruelle--Pollicott} resonances of stochastic systems in reduced state space. {Part I}: Theory},
  journal={Journal of Statistical Physics},
  volume={179}, pages={1366--1402}, year={2020}
}

@article{lin2021data,
  author={Lin, Kevin K. and Lu, Fei},
  title={Data-driven model reduction, {Wiener} projections, and the {Koopman--Mori--Zwanzig} formalism},
  journal={Journal of Computational Physics},
  volume={424}, pages={109864}, year={2021}
}

@article{li2015incorporation,
  author={Li, Zhen and Bian, Xin and Li, Xiantao and Karniadakis, George Em},
  title={Incorporation of memory effects in coarse-grained modeling via the {Mori--Zwanzig} formalism},
  journal={The Journal of Chemical Physics},
  volume={143}, number={24}, pages={243128}, year={2015}
}

@article{hansen2003jensen,
  author={Hansen, Frank and Pedersen, Gert K.},
  title={{Jensen}'s operator inequality},
  journal={Bulletin of the London Mathematical Society},
  volume={35},
  pages={553--564},
  year={2003},
  doi={10.1112/S0024609303002200}
}

@article{chen2024learning,
  author={Chen, Yuan and Xiu, Dongbin},
  title={Learning stochastic dynamical system via flow map operator},
  journal={Journal of Computational Physics},
  volume={508},
  pages={112984},
  year={2024},
  doi={10.1016/j.jcp.2024.112984}
}

\end{document}